%% file: ArtBMS1.tex
\documentclass[referee]{svjour3}

\usepackage{graphicx}
\usepackage{bm}

\input{macroBMS1.tex}

\bd

\title{A discrete time neural network model with spiking neurons}
\subtitle{ Rigorous results on the spontaneous dynamics.}
\titlerunning{Discrete time spiking neurons.}

\author{
B. Cessac
}

\institute{
\at INRIA, 2004 Route des Lucioles, 06902 Sophia-Antipolis, France.
\at INLN, 1361, Route des Lucioles, 06560 Valbonne, France.
\at Universit\'e de Nice, Parc Valrose, 06000 Nice, France.
}

\date{\today}

\maketitle

\input{abstract}

\pagebreak

\input{Intro}

\su{General context.}\label{SDef} 

\input{Mathdef} 
\input{DescBMS}

\su{Preliminary results.}\label{SMF}

\input{SPhase}
\input{SPart}
\input{SPropF}
\input{SSing}

\su{Asymptotic dynamics.}\label{SAS}

We now focus on the asymptotic dynamics of (\ref{DNN}).
\input{SAtt}
\input{SWs}
\input{SMarkov}

\input{SGhost}

\input{SStructA}

\input{SGenA}

\su{Coding dynamics with spiking sequences.} 
\input{SCoding}

\input{SGlob}

\su{Discussion} \label{SDisc} 

\input{SRS}

\input{SNoise}
\input{SIt}

\input{SHebb}
\input{SGibbs}
\input{Sdt}

\input{Ack.tex}

\input{biblio.tex}
\ed

%% file: macroBMS1.tex

\newcommand{\bth}{\begin{theorem}}
\newcommand{\enth}{\end{theorem}}

\newcommand{\bp}{\begin{proposition}}
\newcommand{\ep}{\end{proposition}}
\newcounter{Ccor}[section]
\newtheorem{corrollary}{Corrollary}[Ccor]
\newcommand{\bcor}{\begin{corrollary}}
\newcommand{\ecor}{\end{corrollary}}
\newcommand{\blem}{\begin{lemma}}
\newcommand{\elem}{\end{lemma}}
\newcommand{\bpb}{\begin{question}}
\newcommand{\epb}{\end{question}}

\newcommand{\bdf}{\begin{definition}}
\newcommand{\edf}{\end{definition}}
\newcommand{\bpr} {\begin{proof}}
\newcommand{\epr}{\qed \end{proof}}

\newcommand\bbbr{{\sf I\!R}}

\newcommand\bbbz{{\sf I\!Z}}

\newcommand{\bd}{\begin{document}}
\newcommand{\ed}{\end{document}}
\newcommand{\beq}{\begin{equation}}
\newcommand{\eeq}{\end{equation}}
\newcommand{\bef}{\begin{figure}}
\newcommand{\enf}{\end{figure}}
\newcommand{\bea}{\begin{eqnarray}}
\newcommand{\eea}{\end{eqnarray}}
\newcommand{\baR}{\begin{array}}
\newcommand{\eaR}{\end{array}}
\newcommand{\bc}{\begin{center}}
\newcommand{\ec}{\end{center}}
\newcommand{\ben}{\begin{enumerate}}
\newcommand{\een}{\end{enumerate}}
\newcommand{\bit}{\begin{itemize}}
\newcommand{\eit}{\end{itemize}}
\newcommand{\su}{\section}
\newcommand{\ssu}{\subsection}
\newcommand{\sssu}{\subsubsection}
\newcommand{\nid}{\noindent}
\newcommand{\nnb}{\nonumber}

\newcommand\be{{\bf e}}
\newcommand\bv{{\bf v}}
\newcommand\h{{\bf h}}
\newcommand\bEta{\mbox{{\boldmath $\eta$}}}
\newcommand\cEta{\cC(\bEta)}
\newcommand\bEtap{\mbox{{\boldmath $\eta'$}}}
\newcommand\cEtap{\cC(\bEtap)}
\newcommand\betat{\tilde{\bEta}}
\newcommand\betatp{\tilde{\bEtap}}
\newcommand\betatpl{\tilde{\bEta}^+}
\newcommand\bbeta{\mbox{{\boldmath $\beta$}}}
\newcommand\mub{\mu_{\bbeta}}
\newcommand\mus{\mu^\ast}
\newcommand\Ps{P^\ast}
\newcommand\Fb{\cF_{\bbeta}}
\newcommand\B{{\bf B}}
\newcommand\V{{\bf V}}
\newcommand\Vp{{\bf V'}}
\newcommand\W{{\bf W}}
\newcommand\Z{{\bf Z}}
\newcommand\F{{\bf F}}
\newcommand\Fe{\F_{\bEta}}
\newcommand\Fei{F_{\bEta,i}}
\newcommand\G{{\bf G}}
\newcommand\I{{\bf I}}
\newcommand\bz{{\bf 0}}
\newcommand\bzs{\left\{\bz\right\}}
\renewcommand\S{{\bf S}}
\newcommand\Td{T\left(\dAS\right)}
\newcommand\X{{\bf X}}
\newcommand\x{{\bf x}}
\newcommand\y{{\bf y}}
\newcommand\xt{\X_t(\V)}
\newcommand\cktv{\chi_{k,t}(\V)}
\newcommand\cA{{\cal A}}
\newcommand\cB{{\cal B}}
\newcommand\cH{{\cal H}}
\newcommand\cN{{\cal N}}
\newcommand\cO{{\cal O}}
\newcommand\bme{\cB_\cM(\epsilon)}
\newcommand\Bev{\cB(\V,\epsilon)}
\newcommand\coMe{\stackrel{\circ}{\cMe}}
\newcommand\beV{\cB(\V,\epsilon)}
\newcommand\bmd{\cB_\cM(\delta)}
\newcommand\bmdp{\cB_\cM(\delta')}
\newcommand\cC{{\cal C}}
\newcommand\cD{{\cal D}}
\newcommand\cJ{{\cal J}}
\newcommand\cR{{\cal R}}
\newcommand\De{\cD(\bEta)}
\newcommand\bDe{\bar{\cD}(\bEta)}
\newcommand\Dep{\cD(\bEtap)}
\newcommand\bDep{\bar{\cD}(\bEtap)}
\newcommand\cG{{\cal G}}
\newcommand\GW{\cG_{(\cW,\Ie)}}
\newcommand\IW{\cI_{(\cW,\Ie)}}
\newcommand\cF{{\cal F}}
\newcommand\bcF{\bar{\cF}}
\newcommand\oGF{{\omega(\GF)}}
\newcommand\oGFTe{{\omega(\GFTe)}}
\newcommand\oF{{\omega_{\cF}(\GF)}}
\newcommand\oFTe{{\omega_{\cF}(\GFTe)}}
\newcommand\Bie{\cB_{\infty,\epsilon}}
\newcommand\BTe{\cB_{T,\epsilon}}
\newcommand\BTep{\cB_{{T+1},\epsilon}}
\newcommand\BTeu{\cB_{T,\epsilon}^{(1)}}
\newcommand\BTed{\cB_{T,\epsilon}^{(2)}}
\newcommand\DFFTV{\left. D\FTF\right|_{\V}}
\newcommand\FF{{\F_{\cF}}}
\newcommand\FTF{{\FT_{\cF}}}
\newcommand\FFb{{\F_{\bcF}}}
\newcommand\VF{{\V_{\cF}}}
\newcommand\VFb{{\V_{\bcF}}}
\newcommand\ft{\cF_t(\V)}
\newcommand\ftmu{\cF_{t-1}(\V)}
\newcommand\ftmk{\cF_{t-k}(\V)}
\newcommand\ftmkpu{\cF_{t-k+1}(\V)}
\newcommand\bft{\bar{\cF}_t(\V)}
\newcommand\bftmu{\bar{\cF}_{t-1}(\V)}
\newcommand\bftmk{\bar{\cF}_{t-k}(\V)}
\newcommand\bftmkpu{\bar{\cF}_{t-k+1}(\V)}
\newcommand\cM{{\cal M}}
\newcommand\cMz{{\cal M_{\bf{0}}}}
\newcommand\cMu{{\cal M_{\bf{1}}}}
\newcommand\cMe{{\cM_{\bEta}}}
\newcommand\cMeu{{\cM_{{\bEta}_1}}}
\newcommand\cMet{{\cM_{{\bEta}_t}}}
\newcommand\cMeT{{\cM_{{\bEta}_T}}}
\newcommand\cMep{{\cM_{\bEtap}}}
\newcommand\cU{{\cal U}}
\newcommand\cS{{\cal S}}
\newcommand\cI{{\cal I}}
\newcommand\cP{{\cal P}}
\newcommand\cQ{{\cal Q}}
\newcommand\cV{{\cal V}}
\newcommand\cW{{\cal W}}
\newcommand\Vpi{V^+_i}
\newcommand\vp{\cV^+}

\newcommand\Is{\I^s}
\newcommand\IsP{\I^{s\ast}}
\newcommand\Ise{\I^s(\bEta)}
\newcommand\IE{\I(\bEta)}
\newcommand\Isv{\I^s(\V)}
\newcommand\Isi{I^s_i}
\newcommand\Isj{I^s_j}
\newcommand\Iie{I_i(\bEta)}
\newcommand\Iiep{I_i(\bEta')}
\newcommand\Ije{I_j(\bEta)}
\newcommand\Ijep{I_j(\bEta')}
\newcommand\Isie{\Isi(\bEta)}
\newcommand\Isje{\Isj(\bEta)}
\newcommand\Isjep{\Isj(\bEtap)}
\newcommand\Isiv{\Isi(\V)}
\newcommand\Isvt{I^s_i(\V(t))}
\newcommand\Ie{{\bf I^{ext}}}
\newcommand\Iei{I^{ext}_i}
\newcommand\Iej{I^{ext}_j}
\newcommand\tk{\tau^{(k)}}
\newcommand\tkp{\tau^{(k+1)}}
\newcommand\toi{\tau^{(1)}_i(\V)}
\newcommand\toik{\tau^{(k)}_i(\V)}
\newcommand\toike{\tau^{(k)}_i(\betat)}
\newcommand\toikpe{\tau^{(k+1)}_i(\V)}
\newcommand\toikep{\tau^{(k+1)}_i(\betat)}
\newcommand\dik{\delta_i^{(k)}(\V)}
\newcommand\PrW{\cP_\cW}
\newcommand\PrI{\cP_\Ie}
\newcommand\PrWi{\cP_{\cW,\Ie}}
\newcommand\PW{\cP_{(\cW,\Ie)}}
\newcommand\RW{\cR_{(\cW,\Ie)}}
\newcommand\cPT{\cP^{(T)}}
\newcommand\cMT{\cM_{{\bEta_0} \dots {\bEta_T}}}
\newcommand\FT{\F^T}
\newcommand\Ft{\F^t}
\newcommand\Ftp{\F^{t+1}}
\newcommand\Fttk{\F^{t_k}}
\newcommand\Fmu{\F^{-1}}
\newcommand\Fmd{\F^{-2}}
\newcommand\FmT{\F^{-T}}
\newcommand\dAS{d(\oM,\cS)}
\newcommand\dSV{d(\tVp,\cS)}
\newcommand\trV{\left\{\V(t)\right\}_{t \geq 0}}
\newcommand\tV{\tilde{\V}}
\newcommand\tVp{\tilde{\V}^+}
\newcommand\tpi{t^{pre}_i}
\newcommand\tpj{t^{post}_j}
\newcommand\mW{m_{\cW}}
\newcommand\sW{\sigma_{\cW}}
\newcommand\mI{m_{\Ie}}
\newcommand\sI{\sigma_{\Ie}}
\newcommand\sw{\sigma_{(\cW,\Ie)}}
\newcommand\sB{\sigma_B^2}
\newcommand\sd{\sigma^2}
\newcommand\Cst{[\betat]_{s,t}}
\newcommand\defCst{\left\{\betat | \bEta_s=\as, \dots, \bEta_t = \at    \right\}}
\newcommand\as{\mbox{{\boldmath{$\alpha$}}}_s}
\newcommand\at{\mbox{{\boldmath{$\alpha$}}}_t}
\newcommand\moytpsi{\bar{\psi}_\cW}
\newcommand\sn{\left\{1 \dots N \right\}}
\newcommand\GF{\Gamma_\cF}
\newcommand\GFTe{\Gamma_{\cF,T,\epsilon}}
\newcommand\GFTep{\Gamma_{\cF,T,\epsilon'}}
\newcommand\GFTpe{\Gamma_{\cF,{T+1},\epsilon}}
\newcommand\PF{\Pi_{\cF}}
\newcommand\PbF{\Pi_{\bcF}}
\newcommand\oM{\Omega}

\newcommand\wsg{\cW^s(\V)}
\newcommand\wsl{\cW^s_{loc}(\V)}
\newcommand\wslf{\cW^s_{loc}(\F(\V))}
\newcommand\wse{\cW^s_{\epsilon}(\V)}
\newcommand\wset{\cW^s_{\epsilon}(\V(t))}
\newcommand\tjk{\tau_j^{(k)}}
\newcommand\tik{\tau_i^{(k)}}
\newcommand\tiu{\tau_i^{(1)}}
\newcommand\tikp{\tau_i^{(k+1)}}
\newcommand\Vr{V_{reset}}
\newcommand\Vm{V_{min}}
\newcommand\VM{V_{max}}
\newcommand\Vs{V^{+}}
\newcommand\Vsi{V_i^{+}}
\newcommand\SLp{\Sigma_\Lambda^+}
\newcommand\SL{\Sigma_\Lambda}
\newcommand\SW{\Sigma_{(\cW,\Ie)}}
\newcommand\SWp{\Sigma_{(\cW,\Ie)}^+}
\newcommand{\D}{\displaystyle}
\newcommand{\deq}{\stackrel {\rm def}{=}}
\newcommand{\peq}{\stackrel {\rm \mu p.s.}{=}}

%% file: abstract.tex
\begin{abstract} We derive rigorous results describing the asymptotic dynamics
of a discrete time model of spiking neurons introduced in \cite{BMS}. 
Using symbolic dynamic techniques we show how the dynamics of membrane potential
has a one to one correspondence with sequences of spikes patterns (``raster plots'').
Moreover, though the dynamics is generically periodic, it has a weak form of initial conditions
sensitivity due to the presence of a sharp threshold in the model definition.
 As a consequence, the model exhibits a dynamical regime indistinguishable from chaos
 in numerical experiments.


\end{abstract}

%% file: Intro.tex
\label{Intro}

The description of neuron dynamics can use two distinct representations.
On the one hand, the membrane potential is the  physical variable describing the state
of the neuron and its evolution is ruled by fundamental laws
of physics. On the other hand, a neuron is an excitable medium and its activity
is manifested by emission of action potential or ``spikes'': individual spikes, bursts, spikes trains etc...
The first representation constitutes the basis of  almost all neuron models, and
the Hodgkin-Huxley equations are, from this point of view,
 certainly one of the most achieved mathematical representation of the neuron \cite{HH}.
However, neurons communicate by emission of spikes, and it is likely
that the information is encoded in the
neural code, that is, the sequences of spikes exchanged by the neurons
and their firing times.
Since the spikes emission results from the dynamics of membrane potentials,
the information contained in spikes trains is certainly also
contained in membrane potential dynamics. But switching
from membrane potentials to spikes dynamics allows one to focus on information processing
aspects \cite{Gerstner}.
However, this change of description  is far from being evident, even
when using simple neuron models (see \cite{MB} for a review). Modeling a spike by a certain shape
(Dirac peaks or more complex forms), with a certain refractory period,
etc .. which information have we captured and what have we lost ?
These questions are certainly too complex to be answered in a general setting
(for a remarkable description of spikes dynamics and coding see \cite{Spikes}).

Instead, it can be useful to focus on simplified models of  neural networks,
where the correspondence between the membrane potential dynamics  and spiking
sequences can be written explicitly. This is one of the goals of  the  present work.
 We consider a simple
model of spiking neuron, derived from the leaky integrate and fire model \cite{Gerstner}, but
where the time is discretised. To be the best of our knowledge, this model
has been first introduced by G. Beslon, O. Mazet and H. Soula \cite{HS},\cite{BMS}, and we shall call
it ``the BMS model''. 
Certainly, the simplifications involved, especially the time discretisation,
raise delicate problems concerning  biological interpretations, compared
to more elaborated models or to biological neurons \cite{CV} (see the discussion section).
 But the main interest of the model is its simplicity
and the fact that, as shown in the present paper, one can establish an explicit one-to-one correspondence between the membrane potential
dynamics and the dynamics of spikes. Thus, no information is lost when switching from
one description to the other, even when the spiking sequences have a  complex structure.
 Moreover, this correspondence opens up the possibility of using 
tools from dynamical systems theory, ergodic theory, and statistical physics to address
questions  such as:

\bit
\item How to measure the information content of a spiking sequence ?
\item What is the effect of synaptic plasticity (Long Term Depression, Long Term Potentiation, 
Spike Time Dependent Plasticity, Hebbian learning) 
on the spiking sequences displayed by the neural network ?
\item What is the relation between a presented input and the resulting spiking sequence, before and after learning.
\item What is the effect of stochastic perturbations ? Can we relate the dynamics of the discrete 
time BMS model with noise to previous studies 
on continuous time Integrate and Fire neural networks perturbed by a Brownian noise (e.g. \cite{BrunHak},\cite{RBH}) ?
\eit

This paper is the first one of a series trying to address some of these questions
in the context of BMS model. The goal the present article, is to pose the mathematical framework
used for subsequent developments. 
In section \ref{SMF} we present the BMS model and provide elementary mathematical
results on the system dynamics. We show that the presence
of a sharp threshold for the model definition of neuron firing
 induces singularities responsible for a weak form of
initial conditions sensitivity. This effect is different from the usual notion
of chaos since it arises punctually, whenever a trajectory intersects a zero Lebesgue measure set, called the
 singularity set.
Similar effects are encountered in billiards \cite{Chernov} or in Self-Organized Criticality
\cite{BCK1},\cite{BCK2},\cite{BCKM}. Applying methods from dynamical systems theory we 
 derive rigorous results describing the asymptotic dynamics in section \ref{SAS}.
 Although we show that the dynamics is generically periodic,
the presence of a singularity set has strong effects.
 In particular the number of periodic orbits and the transients
growth exponentially as the distance between the attractor and the singularity set tends to zero.
This   has  a strong impact on the numerics and there is a dynamical regime
numerically indistinguishable from chaos. Moreover, these effects become prominent when perturbing the dynamics or when
the infinite size limit is considered. In this context 
we discuss the existence of a Markov partition allowing to encode symbolically the dynamics
with ``spike trains''. In section \ref{SCoding} we indeed show that there is a one to one correspondence
between the membrane potential dynamics and the sequences of spiking patterns (``raster plots'').
This  opens up the possibility to use methods from ergodic theory and statistical mechanics (thermodynamic
formalism) to  analyse spiking sequences. This aspect will be the central topic of another paper.
As an example, we briefly analyze the case of random
synapses and inputs on the dynamics and compare our analysis to the results obtained
by BMS in  \cite{BMS},\cite{HS}. We exhibit numerically a sharp transition between a neural death regime
where all neurons are asymptotically silent, and a phase with long transient having the appearance
of a chaotic dynamics. This transition occurs for example when the variance of the synaptic
weights increases. A further increase leads to a periodic dynamics with small period.
In the discussion section we briefly comment some extensions (effect of Brownian noise, use of Gibbs
measure to characterize the statistics of spikes) that will be developed in forthcoming papers.\\

\textbf{Warning}
This paper is essentially mathematically oriented (as the title suggests), although some extensive parts
are devoted to the interpretation and consequences of mathematical results for neural networks. Though the proof of theorems
and the technical parts can be skipped, the non mathematician reader interested in computational neurosciences,
may nevertheless have difficulties to find what he gains from this study. Let us briefly comment this point. There is still
a huge distance between the complexity of the numerous models of neurons or neural networks, and the mathematical
analysis of their dynamics, though a couple of remarkable results have been obtained within the 50 past years (see
e.g. \cite{CS} and references therein). This has several consequences and drawbacks. There is a constant
temptation to simplify again and again the canonical equations for the neuron dynamics (e.g. Hodgkin-Huxley
equations) to obtain apparently tractable models. A typical example concerns integrate and fire (IF) models. 
The introduction of  sharp threshold and  instantaneous reset gives a rather simple formulation of neuron activity,
and, at the level of an isolated neuron, a couple of important quantities such 
as the next time of firing can be computed exactly. The IF structure can be extended to conductance
based models \cite{RD,CV} closer to biological neurons.
 However,  there are quite a few rigorous results dealing with the dynamics
of IF models at the \textit{network} level. The present paper provides an example of an IF Neural Network
analysed in a global and rigorous manner. 

The lack of mathematical results concerning
the dynamics of neural networks has other consequences. There is an extensive use of numerical
simulations, which is fine. But the present paper shows the limits of numerics in a model
where ``neurons'' have a rather simple structure. What is for more elaborated models ?
It also warns the reader against the uncontrolled use of  terminologies such as ``chaos, edge of chaos, complexity''. 
In this paper, mathematics allows us to precisely define and analyse mechanisms generating
 initial conditions sensitivity, which are basically presents in all IF neural networks, since they are due to the
sharp threshold. We also give a precise meaning to the ``edge of chaos'' and actually give
a way to locate it.
We evidence mechanisms, such as the first firing of a neuron after an arbitrary large time, which can basically
exist in real neural networks, and raise huge difficulties when willing to decide, experimentally
or numerically, what is the nature of dynamics. Again, what happens for more elaborated models ?
This work is a first step in providing a mathematical setting allowing to handle these questions for
more elaborated IF neural networks models \cite{CV}.

%% file: Mathdef.tex
\ssu{Model definition.} \label{MathDef}

Fix $N > 0$ a positive integer called ``the dimension of the neural network'' (the number
of neurons). 
Let $\cW$ be an $N \times N$ matrix, called ``the matrix of synaptic weights'', 
with entries $W_{ij}$.  It defines
an oriented and signed graph, called ``the neural network associated to $\cW$'',  with  vertices $i=1 \dots N$ called the ``neurons''.
There is  oriented edge $j \to i$ whenever $W_{ij} \neq 0$.  $W_{ij}$ is called ``the synaptic weight from neuron $j$ to neuron $i$''.
The synaptic weight is called ``excitatory'' if $W_{ij}>0$ and ``inhibitory'' if $W_{ij}<0$.

Each vertex (neuron) $i$ is characterized by a real variable $V_i$ called the ``membrane potential of neuron $i$''.
Fix a positive real number $\theta>0$ called the ``firing threshold''. Let $Z$ be the function
$Z(x)=\chi(x \geq \theta)$ where $\chi$ is the indicatrix function. Namely,
$Z(x)=1$ whenever $x \geq \theta$ and $Z(x)=0$ otherwise. $Z(V_i)$ is called the ``firing state of neuron $i$''.
When $Z(V_i)=1$ one says that neuron $i$  ``fires'' and when $Z(V_i)=0$ neuron $i$
is ``quiescent''. Finally, fix $\gamma \in [0,1[$,
called the ``leak rate''. The discrete time and synchronous dynamics of the BMS model is given
by:

\beq \label{DNN}
\V(t+1)=\F(\V(t)),
\eeq

\nid where $\V=\left\{V_i\right\}_{i=1}^N$ is the vector of membrane potentials  and
$\F=\left\{F_i\right\}$ with:

\beq\label{Fi}
F_i(\V)=\gamma V_i \left(1 - Z[V_i] \right)+ \sum_{j=1}^N W_{ij}Z[V_j]+ \Iei; \qquad i=1 \dots N.
\eeq

\nid The variable $\Iei$ is called ``the external current\footnote{\label{notecurrent}From a strict
point of view, this is rather a potential. Indeed, this term is divided by a capacity
$C$ that we have set equal to $1$ (see section \ref{DescBMS} for an interpretation 
of equation (\ref{DNN})). We shall not use this distinction in the present paper.} applied to neuron $i$''.
We shall assume in  this paper that this current does not depend on time
(see however the discussion section from an extension of the
present results to time dependent external currents).
The dynamical system  (\ref{DNN}) is then autonomous.

In the following we shall use the quantity

\beq\label{IsynZ}
\Isiv=\sum_{j=1}^N W_{ij}Z[V_j].
\eeq

\nid  called the ``synaptic current'' received by neuron $i$. The ``total current'' is :

\beq \label{Itot}
I_i(\V)=\Isiv+\Iei
\eeq

Define the firing times of neuron $i$, for the trajectory\footnote{Note that, since the dynamics
is deterministic, it
is equivalent to fix the forward trajectory  or the initial condition $\V\equiv \V(0)$. }  $\V$, by:

\beq\label{Tdech}
\tik(\V)=\inf\left\{t  \ | t > \tau_i^{(k-1)}(\V), \ V_i(t) \geq \theta  \right\}
\eeq  

\nid where $\tau_i^{0}=-\infty$. 

%% file: DescBMS.tex
\ssu{Interpretation of BMS model as a Neural Network.}\label{DescBMS}

The BMS model is based on the evolution equation for the leaky 
integrate and fire neuron \cite{Gerstner} :

\beq\label{IFcont}
\frac{dV_i}{dt}=-\frac{V_i}{\tau} + \frac{I_i(t)}{C}
\eeq

\nid where 
$\tau=RC$ is the integration time scale, with $R$, the membrane resistance and  $C$
the electric capacitance of the membrane. $I_i(t)$ is the  synaptic
 current (spikes emitted by other neurons and transmitted to neuron $i$
via the synapses $j \to i$) and an external current. The equation (\ref{IFcont}) holds whenever the membrane
potential is smaller than a threshold $\theta$, usually depending on time
(to account for characteristics such as  refractory period of the neuron). 
When the membrane potential exceeds the threshold value, the neuron ``fires'' (emission
of an action potential or ``spike''). The spike shape depends on the model.
In the present case, the membrane potential is reset instantaneously to a value $\Vr$,  corresponding to the value of the membrane
potential when the neuron is at rest. More elaborated models can be proposed accounting 
for refractory period, spikes shapes, etc ... \cite{Gerstner}  \\

A formal time discretization of (\ref{IFcont}) (say with an Euler scheme) gives:

\beq 
V_i(t+dt)=V_i(t)\left(1-\frac{dt}{\tau} \right) + \frac{I_i(t)}{C}dt
\eeq

Setting $dt=1$ \footnote{This can be interpreted as choosing the sampling time scale $dt$ smaller
than all characteristic time scales in the model, with similar effects of refractoriness and synchronization. However, this requires a more complete discussion, done in
a separate paper  \cite{CV}. See also section \ref{Sdt}.} and  $\gamma= 1-\frac{1}{\tau}$, we obtain.

\beq \label{IFdisc}
V_i(t+1)=\gamma V_i(t) + \frac{I_i(t)}{C}
\eeq

This discretization imposes that $\tau \geq 1$ in
(\ref{IFcont}), thus $\gamma \in [0,1[$. This equation holds whenever $V_i(t) < \theta$. 
As discussed in e.g. \cite{Iz} it provides a rough but realistic approximation of biological
neurons behaviours.   Note that in biological neurons, 
a spike duration is not negligible but has a finite duration (of order $1$ ms).

The firing of  neuron $i$ is characterized by: 

$$V_i(\tik) \geq \theta$$

and:

\beq\label{Reset}
V_i(\tau^{(k)}_i+1)=\Vr+I_i(\tk)
\eeq 

\nid where, from now on, we shall consider that $C=1$ and that $\Vr$, the reset
potential, is equal to $0$.
Introducing the function $Z$ allows us to write the neuron evolution before
and after firing in a unique equation (\ref{Fi}). Moreover, this apparently naive token provides useful insights
in terms  of symbolic dynamics and interpretation of neural coding.

Note that the firing is not instantaneous. The membrane
potential is maintained at a value  $\theta$ during the time interval $[\tau^{(k)}_i,\tau^{(k)}_i+1[$.
Note also that 
 \textit{simultaneous} firing of several neurons  can occur. 
Moreover, a localized excitation may induce a chain reaction
where  $n_1$ neurons fire at the next time, inducing the firing
of $n_2$ neurons, etc $\dots$. Thus, a localized input may generate
a network reaction on an arbitrary large space scale, in a relatively 
short time scale. The evolution of this propagation phenomenon depends
 on the synaptic weights 
and on the membrane potential values of the nodes involved in the 
chain reaction.
 This effect, reminiscent of the ``avalanches'' observed
in the context of self-organized criticality \cite{Bak}, may have
an interesting incidence in the neural network (\ref{DNN}).

%% file: SPhase.tex
\ssu{Phase space $\cM$.}\label{SPhase}

Since $\gamma < 1$ one can restrict the phase space of (\ref{DNN})
to a compact set\footnote{Note that in the original version of BMS, $V_i  \geq 0$.} 
 $\cM=[\Vm,\VM]^N$  such that $\F(\cM) \subset \cM$ where:

\beq \label{Vmin}
V_{min}=\min(0,\frac{1}{1-\gamma}\left[\min_{i=1 \dots N} \sum_{j | W_{ij}<0} W_{ij}+\Iei \right]),
\eeq

\nid and:

\beq\label{Vmax}
\VM = \max(0,\frac{1}{1-\gamma}\left[ \max_{i=1 \dots N}
 \sum_{j | W_{ij}>0}W_{ij} + \Iei \right]),
\eeq

\nid where we use the convention $\sum_{j \in \emptyset} W_{ij}=0$.
Therefore, $ \sum_{j | W_{ij}<0} W_{ij}=0$ (resp. $ \sum_{j | W_{ij}>0}W_{ij}=0$) if
all weights are positive (resp. negative) and $ \sum_{j | W_{ij}<0} W_{ij} \leq 0$
(resp. $  \sum_{j | W_{ij}>0}W_{ij} \geq 0)$.

This results is  easy to show. Indeed, assume that  for all neurons, $\Vm  \leq V_i \leq \VM$. Then, the membrane
potential of neuron  $i$ at the next iteration is 

$$V'_i = \gamma V_i(1-Z(V_i)) + \sum_{j}W_{ij}Z(V_j) + \Iei.$$

Therefore, 

$$ \gamma \Vm(1-Z(V_i)) + \sum_{j | W_{ij}<0} W_{ij} + \Iei
\leq
V'_i 
\leq \gamma \VM(1-Z(V_i)) + \sum_{j | W_{ij} > 0} W_{ij} + \Iei.$$

If $\Vm<0$ then,

$$\Vm = \gamma \Vm + \D{\min_{i=1 \dots N}}
 \left[\sum_{j | W_{ij}<0} W_{ij} + \Iei\right]
\leq \gamma \Vm(1-Z(V_i)) + \sum_{j | W_{ij}<0} W_{ij} + \Iei \leq V'_i,$$

\nid and if $\Vm=0$, then necessarily 
$\D{\min_{i=1 \dots N}}  \left[\sum_{j | W_{ij}<0} W_{ij} + \Iei\right] \geq 0$
and $V'_i \geq 0 = \Vm$.

Similarly, if $\VM>0$ then,

$$\gamma \VM(1-Z(V_i)) + \sum_{j | W_{ij} > 0} W_{ij} + \Iei \leq 
\gamma \VM+ \D{\max_{i=1 \dots N}}\left[\sum_{j | W_{ij} > 0} W_{ij} + \Iei\right] =\VM.$$

\nid and if $\VM=0$, then necessarily $ 
\D{\max_{i=1 \dots N}}\left[\sum_{j | W_{ij} > 0} W_{ij} + \Iei\right] \leq 0$
and $V'_i \leq 0 = \VM$.

Note that the similar bounds hold if $\Iei$ depends on time.

%% file: SPart.tex
\ssu{Phase space $\cM$.} \label{SPart} 

For each neuron one can decompose the interval $\cI = [\Vm,\VM]$ into
$\cI_0 \cup \cI_1$ with $\cI_0=[\Vm,\ \theta[$, $\cI_1=[\theta,\VM]$. If
 $V \in \cI_0$ the neuron is quiescent, otherwise it fires. 
This splitting induces a partition $\cP$ of $\cM$, that we call the ``natural partition''.
The elements of $\cP$ have the following form. Call $\Lambda=\left\{0,1 \right\}^N$.
Let $\bEta=\left\{\eta_1, \dots, \eta_N \right\} \in \Lambda$. This is a $N$ dimensional vector 
with binary components $0,1$. We call such a vector a \textit{spiking state}.
Then $\cM = \D{\bigcup_{\bEta \in  \Lambda} \cMe}$ where:

\beq\label{PartM}
\cMe = \left\{\V \in \cM \ | \ V_i \in \cI_{\eta_i}  \right\}
\eeq

Equivalently, $\V \in \cMe \Leftrightarrow Z(V_i)=\eta_i, \ i =1 \dots N$.
Therefore, the partition $\cP$ corresponds to classifying the membrane potential
vectors according to their spiking state. More precisely, call:

\beq\label{De}
\De=\left\{i \in \left\{1 \dots N\right\} \ | \ \eta_i =1 \right\},
\eeq

\nid and  $\bDe$ the complementary set $\left\{i \in   \left\{1 \dots N\right\}\ | \ \eta_i =0 \right\}$.
Then, whatever the membrane potential $\V \in \cMe$ the  neurons whose index $i \in \De$ will fire
at the next iteration while the neurons  whose index $i \in \bDe$ will stay quiescent. 
In particular, the synaptic current (\ref{IsynZ}) is fixed by the domain
$\cMe$ since  :

\beq\label{Isyn}
\Isiv \equiv \Isie=\sum_{j \in \De} W_{ij}
\eeq

\nid whenever  $\V \in \cMe$. In the same way we shall write $I_i(\bEta)=  \Isie + \Iei$.

$\cP$ has a simple product structure.
 Its domains are hypercubes (thus they are convex) where the edges are parallels 
to the directions $\be_i$ (basis vectors of $\bbbr^N$). More precisely,
for each $\bEta \in \left\{0,1\right\}^N$,

\beq\label{PstructM}
\cMe=\prod_{i=1}^N I_{\eta_i},
\eeq

\nid where $\prod$ denotes the Cartesian product.

%% file: SPropF.tex
\ssu{Elementary properties of $\F$.}  \label{SPropF}

Some elementary, but essential properties of $\F$, are summarized
in the following proposition. We use the notation

\beq\label{ceta}
\cEta=\sum_{j=1}^N \eta_j=\#\De,
\eeq

\nid for the cardinality of $\De$.
This is the number of neurons that will fire in the next iteration whenever
 the spiking pattern is $\bEta$. 

\bp \label{PFGen}

Denote by $\Fe$  the restriction of $\F$ to the domain $\cMe$.
 Then whatever $\bEta \in \Lambda$,

\ben
\item $\Fe$ is affine and differentiable in the interior of its domain $\cMe$.
\item $\Fe$  is
a a contraction with coefficient $\gamma(1-\eta_i)$
in direction $i$.

\item Denote by $D\Fe$
the Jacobian matrix of $\Fe$. Then $D\Fe$  has $\cEta$ 
zero eigenvalues and $N - \cEta$ eigenvalues $\gamma$. 
\item  Call  $ \Fei$ the $i$-th component of $\Fe$ then
\beq\label{Pprodstruct}
\F(\cMe)=\Fe\left[\prod_{i=1}^N \cI_{\eta_i}\right]
=\prod_{i=1}^N F_{\bEta,i}(\cI_{\eta_i})
\eeq
\nid where $\Fei(\cI_0)$ is the \underline{interval} $[\gamma \Vm + \sum_{j=1}^N W_{ij}\eta_j + \Iei, \gamma \theta + \sum_{j=1}^N W_{ij}\eta_j + \Iei[$
and $\Fei(\cI_1)$ is the  \underline{point} $\sum_{j=1}^N W_{ij}\eta_j+ \Iei$.
 More precisely, if $\cEta=k$, the image
of $\cMe$ is a $N-k$ dimensional hypercube, with faces parallel to the canonical
basis vectors $\be_i$ for all $i \notin D(\bEta)$ and with a volume
$\gamma^{N-k}\left[\theta-\Vm\right]^{N-k}$.  

\een
\ep

According to  item (1) we call the domains $\cMe$, ``domains of continuity''of $\F$. \\
 
\bpr
By definition,  $\forall \V \in \cMe$, $\F_i(\V) = \gamma V_i(1 -\eta_i) +\sum_{j=1}^N W_{ij}\eta_j + \Iei$.
$\F$ is therefore piecewise affine, with a constant $I_i(\bEta)=\sum_{j \in \De} W_{ij}+\Iei$ fixed by the domain $\cMe$.
Moreover $\Fe$  is differentiable on the interior of each domain $\cMe$, with:

\beq \label{DF}
\frac{\partial \Fei}{ \partial V_j} = \gamma\delta_{ij}[1-\eta_i].
\eeq

The corresponding Jacobian matrix is thus  diagonal,
constant in the domain $\cMe$, and its eigenvalues are $\gamma[1-\eta_i]$. Each eigenvalue 
is therefore $0$ if $\eta_i=1$ (neuron $i$
fires) and $\gamma$ if $\eta_i=0$ (neuron $i$ is quiescent). Thus, since $\gamma < 1$,
$\Fe$ is a contraction in each direction $i$. 
Once $\cMe$ has been fixed,
the image of each coordinate $V_i$ is only a function of $V_i$. Thus, 
if $V \in \cMe = \prod_{i=1}^N \cI_{\eta_i}$, then
$\Fei(\V)=\Fei(V_i)$ and
$\Fe$ maps the hypercube $\cMe=\prod_{i=1}^N \cI_{\eta_i}$ 
onto the hypercube $\prod_{i=1}^N \Fei(\cI_{\eta_i})$. The segments $\cI_{\eta_i}$
with $\eta_i=0$ are mapped to parallel segments $[\gamma \Vm + \sum_{j=1}^N W_{ij}\eta_j + \Iei,
 \gamma \theta + \sum_{j=1}^N W_{ij}\eta_j + \Iei[$ while each
segment  $\cI_{\eta_i}$
with $\eta_i=1$ is mapped to a point. Thus, if $\cEta=k$ the image
of $\cMe$ is a $N-k$ dimensional hypercube, with faces parallel to the canonical
basis vectors $\be_i$, where $i \notin D(\bEta)$ and with a volume
$\gamma^{N-k}\left[\theta-\Vm\right]^{N-k}$.  
\epr

Finally, we note the following property. 
The dynamical system (\ref{DNN}) can be defined on $\bbbr^N$ and
the contraction property extends to this space. If one considers
the $\delta$-ball $\bmd = \left\{\V \in \bbbr^N | d(\V,\cM)< \delta  \right\}$ then :

\beq \label{ContBoule}
\F\left[\bmd\right]   \subset \bmd.
\eeq
 The distance $d$  is, for example :

\beq\label{dist}
d(\X,\X')=\max_{i=1 \dots N}|X_i - X'_i|,
\eeq

\nid natural in the present context according to property
\ref{PFGen} (eq. (\ref{Pprodstruct})).

%% file: SSing.tex
\ssu{The singularity set $\cS$.}\label{SSing}

The set 

\beq \label{S} 
\cS=\left\{\V \in \cM, \ | \exists i, \ V_i = \theta  \right\},
\eeq

\nid is called the \textit{singularity set} for the map $\F$. 
$\F$ is discontinuous on $\cS$. This set has a simple structure:
this is a finite union of $N-1$ dimensional hyperplanes corresponding
to faces of the hypercubes $\cMe$. Though $\cS$ is a ``small'' set
both from the topological (non residual set) and metric  (zero Lebesgue measure) point of view,
it has an important effect on the dynamics. 

Indeed, let us consider the trajectory of a point $\V \in \cM$ and perturbations
with an amplitude $< \epsilon$ about $\V$. Equivalently, consider the 
evolution of the $\epsilon$ ball $\Bev$ under $\F$.
If $\Bev \cap \cS = \emptyset$ then by definition $\Bev \subset \coMe$, some $\bEta$, where $\coMe$
is the interior of the domain $\cMe$. Thus, by prop. \ref{PFGen}(2) $\F[\cB(\V,\epsilon)] \subset \cB(\F(\V),\gamma\epsilon)$.
More generally, if the images of $\Bev$ under $\Ft$ never intersect $\cS$, then, at time $t$,
 $\F^t[\cB(\V,\epsilon)] \subset \cB(\F^t(\V),\gamma^t\epsilon)$. Since $\gamma<1$, there is a contraction
of the initial ball, and the perturbed trajectories about $\V$ become asymptotically indistinguishable
from the trajectory of $\V$. (Actually, if all neurons have fired after a finite time $t$ then
all perturbed trajectories collapse onto the trajectory of $\V$ after $t+1$ iterations).

On the opposite, assume that there is a time, $t_0$ such that $\F^{t_0}(\Bev) \cap \cS \neq \emptyset$.
By definition, this means that there exists a subset of neurons $\left\{i_1, \dots, i_k\right\}$ and
  $\V'  \in \Bev$, such that $Z(V_i(t_0))\neq Z(V'_i(t_0))$, $i \in \left\{i_1, \dots, i_k\right\}$.
Then:

$$
\baR{lr}
F_i(\V(t_0))-F_i(\V'(t_0)) =&\\ 
\gamma\left[V_i(t_0)(1-Z(V_i(t_0)))-V'_i(t_0)(1-Z(V'_i(t_0)))\right]+
\sum_{j \in \left\{i_1, \dots, i_k\right\}} 
W_{ij}\left[Z(V_j(t_0))-Z(V'_j(t_0))\right]&
\eaR$$

In this case, the difference between $F_i(\V(t_0))-F_i(\V'(t_0))$ is not proportional to $V_i(t_0) - V'_i(t_0)$
, for $i \in \left\{i_1, \dots, i_k\right\}$. Moreover, this distance is finite while $|V_i(t_0) - V'_i(t_0)|<\epsilon$ can
be arbitrary small. Thus, in this case, the crossing of $\cS$ by the $\epsilon$-ball induces
a strong separation effect reminiscent of initial condition sensitivity in chaotic dynamical
system. But the main difference with chaos is that the present effect occurs only
when the ball crosses the singularity. (Otherwise the ball is contracted). The result is a weak
form of initial condition sensitivity and unpredictability occurring also in billiards \cite{Chernov}
or in models of self-organized criticality \cite{BCK1},\cite{BCK2}.
Therefore, $\cS$ is the only source of complexity of the BMS model, and its existence is due to the strict
threshold in the definition of neuron firing. 

Note that if one replaces the sharp threshold
by a smooth one (this amounts to replacing an Heaviside function by a sigmoid) then
the dynamics become expansive in the region where the slope of the regularized threshold
is larger than $1$. Then, the model exhibits chaos in the usual sense (see e.g. \cite{PD},\cite{JP}).
Thus, in some sense, the present model can be viewed as a limit of a regular neural
network with a sigmoidal transfer function. However, when dealing with asymptotic
dynamic one has to consider two limits ($t \to +\infty$ and slope $\to  +\infty$)
that may not commute.

%% file: SAtt.tex
\ssu{The $\omega$-limit set.}\label{SAtt}

\bdf  (From \cite{KH,GH})
A point $y \in \cM$ is called an $\omega$-limit point
for a point $x \in \cM$ if there exists a sequence of times
$\left\{t_k\right\}_{k=0}^{+\infty}$, such that
$x(t_k) \to y$ as  $t_k \to +\infty$. The $\omega$-limit set of $x$, $\omega(x)$,
is the set of all $\omega$-limit points of $x$. The $\omega$-limit
set of $\cM$, denoted by $\oM$, is the set $\oM = \bigcup_{x \in \cM}\omega(x)$.
\edf

Equivalently, $\oM$ is  the set of accumulation points of $\Ft(\cM)$.
In the present case, since $\cM$ is closed and invariant, we have 
$\oM=\bigcap_{t=0}^\infty \F^t(\cM)$.

 The notion of $\omega$ limit set is less known and used than
the notion of \textit{attractor}. There are several
distinct definition of attractor. For example, according to \cite{KH}:

\bdf
A compact set $\cA \in \cM$ is called an \textit{attractor}
for $\F$ if there exists a neighborhood $\cU$ of $\cA$
and a time $N>0$ such that $\F^N(\cU) \subset \cU$ and 

\beq\label{A}
\cA=\D{\bigcap_{t=0}^\infty \F^t(\cU)}.
\eeq 
\edf

Note that from equation (\ref{ContBoule}) one may choose for $\cU$ any open set  such that:
\beq\label{BA}
\cU \supset \bmd, \ \forall \delta >0.
\eeq

In our case $\cA$ and $\oM$ coincide \textit{whenever $\cA$ is not empty}. However, there
are cases where the attractor is empty while the $\omega$ limit
set is not (see example of Fig. 3.3.1 in \cite{KH}, page 128). 
We shall actually encounter the same situation in section
\ref{SGhost}. For this reason we shall mainly use the notion
of $\omega$-limit set instead of the notion of attractor, though we shall
see that they coincide except for a non generic set of synaptic
weights and external currents.

%% file: SWs.tex
\ssu{Local stable manifolds.}\label{SWs}

The stable manifold of $\V$ is the set:

\beq\label{WS}
\wsg=
\left\{
\V' \ | 
d\left(
\F^t(\V'),\F^t(\V)
\right)
\to 0 \quad t \to +\infty
  \right\}.
 \eeq

The local stable manifold $\wsl$ is the largest connected
component of  $\wsg$ containing $\V$. It obeys:

\beq\label{WSloc}
\F\left[\wsl\right] \subset \wslf.
\eeq

In the present model, if  $\V$ has a local stable manifold
$\wse$ of diameter $\epsilon$ then:

\beq\label{Wsinc}
\F^t\left[\wse\right] \subset \cW^s_{\gamma^t \epsilon}(\F^t(\V)).
\eeq

Thus, a perturbation 
of amplitude $<\epsilon$ is exponentially damped and
the asymptotic dynamics of  any point
belonging to the local stable manifold of $\V$
is indistinguishable from the evolution of $\V$.\\

In  BMS model some point may not have a local stable
manifold, due to the presence of the singularity set. 
Indeed, if a small ball of size $\epsilon$ and center $\V$
intersects $\cS$ it will be cut into several pieces 
strongly separated by the dynamics. If this happens,
$\V$ does not have a local stable manifold of size $\epsilon$.
According
to (\ref{Wsinc}) a point $\V \in \cM$ has a local stable manifold
of diameter $\epsilon$ if :

\beq\label{CBC}
\V \notin \bigcap_{t_0 \geq 0} \bigcup_{t \geq t_0} \F^{-t}(\cU_{\gamma^t \epsilon}(\cS)),
\eeq 

\nid where $\cU_\delta(\cS)=\left\{\V \ | \ d(\V,\cS) < \delta \right\}$ is the
$\delta$-neighborhood of $\cS$. This means that the dynamics contracts the $\epsilon$
ball faster than it approaches the singularity set. A condition like (\ref{CBC}) is useful for
measure-theoretic estimations of the set of points having no stable manifold
via the Borel-Cantelli lemma.

In the present context, a more direct approach consists in computing: 

\beq\label{dSV}
\dSV= \inf_{t \geq  0} \min_{i =1 \dots N} |V_i(t) - \theta|,
\eeq

\nid which measures the ``distance'' between the forward trajectory $\tVp \deq \trV$ of $\V$ and $\cS$.
One has the following:

\bp\label{PWSloc}
If  $\dSV > \epsilon > 0$ then $\V$ has a local stable manifold of diameter
$\epsilon$.
\ep

\bpr
This results directly from proposition \ref{PFGen}. Indeed, if
$\dSV>\epsilon$, the image of the $\epsilon$-ball $\Bev$ under $\F^t$, belong
to a unique continuity domain of $\F$, $\forall t >0$ and $\F$ is contracting
on each domain of continuity.
\epr

In the same way, one defines the distance\footnote{Note that this is not a proper distance,
since one may have $d(A,B)=0$ and $A \neq B$. The fact that
$\dAS=0$ if and only if $\oM \cap \cS \neq \emptyset$ is true only because
both sets are closed. I thank one referee for this remark.} of 
the omega limit set $\oM$ to the singularity set 
(one may also consider the distance to the attracting set whenever
$\cA$ is not empty):

\beq\label{dAS}
\dAS= \inf_{\V \in \oM} \dSV.
\eeq

The distance vanishes if and only if $\oM \cap \cS \neq \emptyset$.
Thus, if $\dAS > \epsilon > 0$ any point of $\oM$ has a local stable manifold.
In this situation, any $\epsilon$- perturbation about  $\V \in \oM$ is asymptotically damped.
Note however that  $ \dAS$ can be positive but arbitrary small (see section \ref{SRS}).

%% file: SMarkov.tex
\ssu{Symbolic coding and Markov partition.} \label{SMarkov}

The  partition $\cP$ provides a natural way for
encoding the dynamics. Indeed, to each forward trajectory 
$\tVp$
one can associate an infinite sequence of spiking patterns 
$\bEta_1, \dots, \bEta_t \dots$ where $\bEta_t= \left\{\eta_{i;t}=Z(V_i(t))\right\}_{i=1}^N$.
This sequence provides exactly the times of firing for each neuron.
It contains thus the ``neural code'' of the BMS model. In fact, this sequence is
exactly what biologists call the ``raster plot'' \cite{Gerstner}. On the other hand,
knowing the spiking sequence and the initial condition $\V\equiv \V(0)$
one can determine $\V(t)$ since:

\beq\label{Vit} 
V_i(t)=\gamma^t \prod_{k=0}^{t-1} \left(1-\eta_{i;k}\right)V_i(0) 
+ \sum_{n=1}^t \gamma^{t-n} 
\prod_{k=n}^{t-1} (1-\eta_{i;k})
I_i(\bEta_{n-1}),
\eeq
 
\nid where $I_i(\bEta_{n-1})=\sum_{j=1}^N W_{ij} \eta_{j;n-1}+\Iei$
and where we used the convention
 $\gamma^{t-n} \prod_{k=n}^{t-1} (1-\eta_{i;k})=1$ if $n=t$.
(Note that the same equation holds if $\Iei$ depends on time).
 
The term $\gamma^t \prod_{k=0}^{t-1} \left(1-\eta_{i;k}\right)V_i(0) $ contains
 the initial condition, but it 
vanishes as soon as $\eta_{i;k}=1$, some $k$, (which means that the neuron
has fired at least once between time $0$ and $t-1$).
If the neuron does not fire then   this term is asymptotically
damped. Thus, one can expect that after a sufficiently long time (of order $\frac{1}{|\log(\gamma)|}$),
 the system ``forgets''
its initial condition. Then, knowing the evolution
of $\V(t)$ should be equivalent to knowing the neural code. However,
this issue requires a deeper inspection  using
symbolic dynamics techniques and we shall see that the situation is a little bit more complex
than expected.\\

For this, one first defines a transition graph $\GW$  from the natural partition $\cP$.
This graph depends on the synaptic weights (matrix $\cW$)
and on the external currents (vector $\Ie$) as well.
The vertices of $\GW$ are the spiking patterns $\bEta \in \Lambda=\left\{0,1\right\}^N$.
Thus, one associates to each spiking pattern $\bEta$ a vertex in $\GW$.  
Let $\bEta,\bEta'$ be two vertices of $\GW$. Then there is an oriented
edge $\bEta \to \bEtap$ whenever  $\F(\cMe)  \cap \cMep \neq \emptyset$.
The transition $\bEta \to \bEtap$ is then called \textit{legal}. 
Equivalently, a  legal transition  satisfies the \textit{compatibility conditions}:

\beq\label{Transitions}
\baR{ccccc}
&(a) \quad & i \in \De \cap \Dep& \Leftrightarrow& \sum_{j \in \De} W_{ij} +  \Iei \geq \theta\\ 
&(b) \quad & i \in \De \cap \bDep& \Leftrightarrow& \sum_{j \in \De} W_{ij} + \Iei < \theta \\ 
&(c) \quad & i \in \bDe \cap \Dep& \Leftrightarrow& \gamma V_i + \sum_{j \in \De} W_{ij} + \Iei \geq \theta \\
&(d) \quad & i \in \bDe \cap \bDep& \Leftrightarrow& \gamma V_i +\sum_{j \in \De} W_{ij} + \Iei < \theta 
\eaR
\eeq

\nid (recall that $\De$ is given by eq. (\ref{De})).
The transition graph depends therefore on the coupling matrix $\cW$ and the external current
$\Ie$. It also depends
on the parameters $\gamma,\theta$ but we shall omit this dependence
in the notation. 
Note that the transitions (a), (b) do not  depend on the membrane potential.
We denote by $ \SWp$ the set of right infinite legal sequences $\betatpl=
\left\{\bEta_1, \dots, \bEta_t \dots  \right\}$ and by
$\SW$ the set of bi-infinite sequences $\betat=
\left\{\dots \bEta_s, \dots, \bEta_{-1} \bEta_0 \bEta_1, \dots, \bEta_t \dots \right\}$.\\
 
This coding is particularly useful if there is a one to one correspondence
(except for a negligible set) between a legal sequence and an orbit
of  (\ref{DNN}). This is not necessarily the case due to the presence
of the singularity set. However one has this correspondence whenever one can
construct a \textit{finite Markov partition} by a suitable
refinement of $\cP$. In the present context where the dynamics
is not expanding and just contracting, a partition $\cQ$
is a Markov partition if its elements satisfy
$\F(\cQ_n) \cap \cQ_{n'} \neq \emptyset \Rightarrow \F(\cQ_n)  \subset  \cQ_{n'}$.
In other words, the image of $\cQ_n$ is included in $\cQ_{n'}$
  whenever the transition $n \to n'$ is legal.  \\

$\cP$ is in general not a Markov partition  (except if $\gamma=0$
and maybe for a non generic set of $W_{ij},\Iei$ values). This is because
the image of a domain  usually intersects several domains.
(In this case the image intersects the singularity set). From the neural networks
point of view this means that it is in general not possible to know what
will be the spiking pattern at time $t+1$ knowing the spiking pattern at time $t$.
There are indeed several possibilities depending on the \textit{membrane potential values}
and not only on the firing state of the neurons.
The question is however: knowing a sufficiently large (but finite) sequence
of spiking patterns is it possible,  under some circumstances, to predict which spiking patterns
will come next ? The answer is yes.

\bth \label{TAfini}
Assume that $\dAS>\epsilon>0$.  Then:

\ben

\item Call $\F^t$ the $t$-th iterate of $\F$. There 
is a finite $T$, depending on $\dAS$,  such that $T \to +\infty$ when $\dAS \to 0$
and such that there exists  a finite Markov partition for $\FT$.

\item $\oM$ is a finite union of stable periodic orbits with a finite period.
These orbits are encoded by a sequence of finite blocs of spiking patterns, each bloc
corresponding to a Markov partition element.   

\een

\enth

\bpr  Fix $T >0$. Consider
the partition  $\cPT$ whose elements have the form:
\beq\label{AtPt}
\cMT=\cM_{\bEta_0} \cap \Fmu\left(\cM_{\bEta_1} \right)
\cap \Fmd\left(\cM_{\bEta_2} \right)
\cap \dots
\cap \FmT\left(\cM_{\bEta_T} \right).
\eeq
 By construction
$\FT$ is continuous  and
thus is a contraction from the interior of each domain $\cMT$ into $\cM_{\bEta_T}$, with
$|\FT(\cMT)| \leq \gamma^T |\cMT|$, where $|\cMT| < |\cM_{\bEta_0}|$ and where $| \ |$ denotes
the diameter.
 Thus there  is a finite

\beq\label{T}
T =\left[\frac{\log(\epsilon)-\log(|\cM_{\bEta_0}|)}{\log(\gamma)}\right]
\geq \frac{\log(\dAS)-\log(|\cM_{\bEta_0}|)}{\log(\gamma)},
\eeq

\nid where $\left[ \ \right]$ is the integer part,  such that 
$\forall \cMT$, $|\FT(\cMT)| \leq \epsilon < \dAS$. Then $\cPT$
has finitely many domains ($2^{NT}$).  
Denote them by $\pi_n, \ n =1  \dots 2^{NT}$. 
Then, $|\FT(\pi_n)| \leq \epsilon, \forall n$.

Since $\FT(\oM \cap \pi_n) \subset \oM \cap \FT(\pi_n)$ the points
belonging to $\oM \cap \pi_n$ are mapped, by $\FT$, into
a subset of $\oM$ of diameter $\leq \epsilon$. 
Since $\dAS > \epsilon > 0$ each point in $\oM$ has a  local
 stable manifold of diameter $\epsilon$. Thus all
 points of $\FT(\oM \cap \pi_n)$ belong to the same stable manifold.
Hence all these points converge to the same orbit in $\oM$ and $\pi_n$
contains at most one point in $\oM$. Since there are finitely many
domains $\pi_n$, $\oM$ is composed by finitely many points and since
the dynamics is deterministic, $\oM$ is a finite union of stable periodic orbits with a  finite period.
  If $\pi_n \cap \oM = \emptyset$ then this domain is, by definition, non recurrent
and it is mapped into a union of domains $\pi_{n_k}$ containing a point of $\oM$.  
For all $\pi_n$ containing a point of $\oM$,
$\FT(\pi_n) \cap  \pi_{n'}  \neq \emptyset \Rightarrow \FT(\pi_n) \subset \pi_{n'}$.
Therefore, $\cPT$ is a Markov partition for the mapping $\FT_\Omega$.
\epr

\textbf{Remarks.} 

\bit

\item{\textbf{Structural stability.}} There is a direct consequence of the previous
theorem. Assume that we make a small perturbation of some $W_{ij}$'s or $\Iei$'s.
This will result in slight change of the domains of continuity of $\cP$ and leads
to a perturbed natural partition $\cP'$. This will
also change the $\omega$-limit set. Call the perturbed $\omega$-limit set $\Omega'$.
If $\dAS>\epsilon>0$ then if the perturbation is small enough such
that, for any orbit in $\Omega$, the perturbed and unperturbed orbit have the same
sequence of spiking patterns, then the set $\Omega$
and $\Omega'$ have the same number of fixed points and their distance remains small
(it vanishes when the amplitude of the perturbation tends to zero). This corresponds
to a structurally stable situation. On the opposite, when increasing continuously
the amplitude of the perturbation, there is a moment where  the perturbed and unperturbed orbit
have a different sequence of spiking patterns. This corresponds to a bifurcation 
in the system and the two $\omega$-limit sets can be drastically different.    

\item{\textbf{Maximal period.}} The number

\beq\label{Td}
T_d=2^{N\frac{\log(\dAS)}{\log(\gamma)}},
\eeq

\nid gives an upper bound for the number of Markov partition
elements, hence for the cardinality of $\Omega$ and for the maximal period. 
It increases \textit{exponentially} with the system size $N$
and with $\log(\gamma)$ and $\log(\dAS)$. 
(Note that this time is useful essentially when $\dAS$ is \textit{small} (and lower than $1$)).
 Hence, even if the dynamics is  periodic it can nevertheless
be quite a bit complex.
\eit

Theorem \ref{TAfini} opens up the possibility
of associating to each orbit in $\oM$ a symbolic orbit
constituted by a finite sequence of spiking patterns, whenever $\dAS > \epsilon >0$.
This result is generalized in the section \ref{SGlob}
and its consequence are discussed.

%% file: SGhost.tex
\ssu{Ghost orbits.} \label{SGhost} 

Before proceeding to the  characterisation of the $\omega$-limit set structure
in the general case,
we have to treat a specific situation, where a neuron takes an arbitrary large time to
fire. This situation may look strange from a practical point of view,
but it has deep implications. Indeed, assume that we are in a situation
where we cannot bound the first time of firing of a neuron.
This means that we can observe the dynamics on arbitrary long times
without being able to predict what will happen later on, because when
this neuron eventually fire, it may drastically change
the evolution. This case is exactly related to the chaotic
or unpredictable regime of BMS model. From a mathematical point of view
it may induce ``bad'' properties such as an 
empty attractor. We shall however see that this situation is non generic.

\bdf\label{DfGhost}
An orbit $\tV$ is a \textit{ghost orbit} if  $\exists i$ such that:

$$(i) \forall t> 0, V_i(t)<\theta$$ 

\nid and :

$$(ii) \limsup_{t \to +\infty} V_i(t)=\theta$$
\edf

\textbf{Examples.}

\ben

\item One neuron ($N=1$), $W_{11}=0$, $\Vr=0$ 
and $I^{ext}_1=\theta(1-\gamma)<\theta$. Take $V_1(0)=0$. 
Then, from eq. (\ref{Vit}),
 $V_1(t)=\sum_{n=1}^t \gamma^{t-n} I^{ext}_1 = 
\theta(1 - \gamma^t)< \theta$
and $\lim_{t \to +\infty} V_1(t)=\theta$. Therefore
the orbit of $0$ is a ghost orbit. If $V_1(0) \geq \theta$
the neuron fires and $V_1(1)=I^{ext}$. Thus this point is mapped into 
$\cM= \left[0,I^{ext}\right]$. If  $0 \leq V_1(0) < \theta$
then, $V_1(t)=\gamma^t V_1(0) +
\theta(1 - \gamma^t)$ 
and the neuron fires after a finite time, but then it is mapped
to $V_1=0$. Thus all points of $\cM= \left[0,I^{ext}\right]$ are eventually 
mapped to $0$ and the orbit of $0$ is a ghost orbit. 
In this case $\oM=\left\{0\right\}$ while $\cA$ is empty (see \cite{KH} page 128 for a similar example).

\item Two neurons  with $W_{22}>\theta; 0<W_{12} \geq (1-\gamma)\theta; \ W_{21} >0$
and where for simplicity we assume that $\Vm=0$ ($W_{11}\geq 0$) and $\Iei=0$.
In this case, if $2$ fires once, it will fire forever. Then the dynamics of $1$ is
$V_1(t+1)=\gamma V_1(t)+ W_{12}$, as long as $V_1(t) < \theta$. Therefore, if
$V_1(0) < \theta$, then $V_1(t+1)=\gamma^{t+1} V_1(0) + W_{12}\frac{1-\gamma^{t+1}}{1-\gamma}$
as long as $V_1(t) < \theta$. The condition $V_1(t)<\theta$ is equivalent to $V_1(0) < f(t)$, with
$f(t)=\frac{\theta}{\gamma^t} + \frac{W_{12}}{1-\gamma}(1-\frac{1}{\gamma^t})$.
This function is strictly decreasing if $W_{12}>(1-\gamma)\theta$ and $f(t) \to -\infty$
as $t \to \infty$.
 Thus, for a fixed $W_{12} > (1-\gamma) \theta$ there is a 
$\tau=\left[\frac{\log(1-\frac{\theta(1-\gamma)}{W_{12}})}{\log\gamma}\right]$
(where $[\ ]$ is the integer part), such that $\forall 0 \leq t < \tau$,
there exists and interval $\cJ_t = [f(t),f(t-1)[ \in [0,\theta]$
such that $\forall V_1(0) \in \cJ_t$, the neuron $1$ will fire for the first time at time $t$.
When $W_{12} \to \theta(1-\gamma)$ from above, $\tau$ diverges and one can find an initial condition
such that the first firing time of $1$ is arbitrary large (transient case). This
 generates a ghost orbit. 
\een 

One may generalize these examples to arbitrary dimensions.
However, the previous examples look where very specific since
we had to adjust the parameters to a precise value, and the ghost
orbit can be easily removed by a slight variation of these parameters.
This suggests us that this situation is non generic. We shall
prove this in section \ref{SStructA}. 

To finish this section let us emphasize that, though ``strict''
 ghost orbits, having the limit $t \to \infty$ in the definition,
are non generic, it may happen that $V_i(t)$ remains below
the threshold during an arbitrary long (but finite) time before firing. Then,
the characterization of the asymptotic dynamics
may be out of numerical or experimental control.

%% file: SStructA.tex
\ssu{Two theorems about the structure of $\oM$.} \label{SStructA}

The condition $\dAS>\epsilon > 0$ excludes situations where some points accumulate 
on the singularity set.
In these situations, the usual behavior is the following. An $\epsilon$-ball containing
a point $\V$ accumulating on $\cS$ will be cut in several pieces when it intersects
the singularity set. Then, each of these pieces may intersects $\cS$ later on, etc...
At each intersection the dynamics generates distinct orbits and  strong separations
 of trajectories. 
It may happen that the proliferation of orbits born from an $\epsilon$-ball
 goes on forever and
 there are examples of such dynamical system having  a positive (topological) entropy even if
 dynamics 
is contracting \cite{Rypdal}. Also,  points accumulating on $\cS$ do not have
a local stable manifold.

In BMS model the situation is however less complex, due to the reset term
$\gamma V_i(1-\eta_i)$. Indeed, consider the image of an $\epsilon$ ball $\beV$ about
some point $\V$.  Assume that the ball intersects several domains of continuity. Then,
the action of $\F$ generates several pieces, as in the usual case.
But, the image of $\beV \cap \cMe$ is a $N-\cEta$ dimensional domain,
 whose projection  in each direction
$i$ such that $\eta_i=1$ is a \underline{point}. Thus, even if $\beV$
intersects the $2^N$ domains of $\cP$, its image will be an union
of $2^N$ pieces  all but one having a dimension $<N$.  
This effect limits the proliferation of orbits and the complexity
of the dynamics and the resulting structure of the $\omega$-limit set is relatively simple, even
if $\dAS=0$ provided 
one imposes some additional assumptions. More precisely, the following holds.

\bth\label{Tomega} 
Assume that $\exists \epsilon >0$ and $ \exists T < \infty$ such that,
$\forall \V \in \cM$, $\forall i \in \left\{1 \dots N \right\}$,
\ben
\item Either $\exists t \leq T$ such that $V_i(t) \geq  \theta$;
\item Or $\exists t_0 \equiv t_0(\V,\epsilon)$
 such that $\forall t \geq t_0$, $V_i(t) < \theta -\epsilon$
\een
Then,  $\oM$ is composed by finitely many
periodic orbits with a finite period. 
\enth 

Note that conditions (1) and (2) are not disjoint. The meaning
of these conditions is the following. We impose that
either a neuron have fired after a finite time 
(uniformly bounded, i.e. independent of $\V$) or, if it does not fire after a certain time
it stays  bounded below the threshold value (it cannot accumulate on $\theta$).
 Under these assumptions the asymptotic dynamics
is periodic and one can predict the evolution after observing the system
on a finite time horizon $T$, whatever the initial condition.
Note however that $T$ can be quite a bit large.\\
  
The proof uses the following lemma.

\blem\label{Lomega}
Fix $\cF$ a subset of $ \left\{1 \dots N \right\}$ and let $\bcF$ be the complementary set of $\cF$. Call
$$\GFTe=
\left\{
\V \in \cM \left| 
\baR{ccc}
&(i) \ \forall i  \in \cF,&  \exists t \leq T, \mbox{such \ that} \ V_i(t) \geq \theta\\
&(ii) \ \forall j \in \bcF,& \exists t_0 \equiv t_0(\V,j) <  \infty ,  \mbox{such \ that} 
\ \forall t > t_0,  V_j(t) < \theta -\epsilon
\eaR
\right.
\right\}
$$
\nid then $\omega(\GFTe)$, the $\omega$-limit set of $\GFTe$,  is composed by finitely many
periodic orbits with a finite period. 
\elem

\bpr of th. \ref{Tomega}

Note that there are finitely many subsets $\cF$ of $ \left\{1 \dots N \right\}$.
Note also that $\GFTe \subset \GFTpe$ and that $\GFTe \subset \GFTep$ whenever $\epsilon' < \epsilon$. 
We have therefore:

$$\cM \subset \bigcup_{\cF} \bigcup_{T>0} \bigcup_{\epsilon >0} \GFTe 
= \bigcup_{\cF}
 \Gamma_{\cF,+\infty,0}.$$

But, under hypothesis (1) and (2) of th. \ref{Tomega}, there exists $\epsilon > 0, T <\infty $
such that $\cM=\bigcup_{\cF} \GFTe$ where the union on $\cF$ is finite.
Since $\F(\cM) \subset \bigcup_{\cF} \F(\GFTe)$, $\oM \subset  \bigcup_{\cF} \omega(\GFTe)$.
Under lemma \ref{Lomega} 
$\oM$ is therefore a subset of a finite union of sets
containing  finitely many periodic orbits with a finite period.
\epr

\bpr of lemma \ref{Lomega}
Call $\PF$ (resp. $\PbF$) the projection onto the subspace generated by the basis vectors $\be_i, \ i \in \cF$
(resp.  $\be_j, \ j \in \bcF$) and set $\VF = \PF\V$ ($\VFb=\PbF \V$), $\FF=\PF\F$ ($\FFb=\PbF \F$).
Since each neuron $j \in \bcF$ is such that:
\beq\label{VjF}
V_j(t) = \sum_{n=0}^{t-t_j-1} \gamma^{n} (\sum_{k} W_{jk}Z[V_k(t-n-1)]+I_j^{ext}) 
< \theta-\epsilon,
\eeq
for $t$ sufficiently large, (larger than the last (finite) firing time $t_j$),
  these neurons do not act on the other neurons 
and their membrane potential
is only a function of the synaptic current generated by the neurons $\in \cF$. 
Thus, the asymptotic dynamics 
is generated by the neurons $i \in \cF$. Namely,
$\forall \V \in \oGFTe$, $\VF(t+1)=\FF[\VF(t)]$ and
$\VFb(t+1)=\FFb[\VF(t)]$. One can therefore focus the analysis of
the $\omega$ limit set to its projection $\oFTe=\PF\oGFTe$ (and infer the dynamics
of the neurons $j \in \bcF$  via eq. (\ref{VjF})). 

Construct now the partition $\cPT$, with convex elements
given by $\cMT=\cM_{\bEta_0} \cap \Fmu\left(\cM_{\bEta_1} \right)
\cap \Fmd\left(\cM_{\bEta_2} \right)
\cap \dots
\cap \FmT\left(\cM_{\bEta_T} \right)$, where $T$ is the same as
in the definition of $\GFTe$. By construction, $\FT$ is continuous on each element of $\cPT$
and fixing $\cMT$ amounts to fix the affinity constant of $\FT$. 
By definition of $T$, $\DFFTV$,    the derivative of $\FTF$ at $\V$, 
has all its eigenvalues equal to $0$ whenever $\V \in \oFTe$ (prop. \ref{PFGen}.3). 
Therefore $\FTF[\cMT \cap \oFTe]$
is a point. Since 
$$\FTF(\cM \cap \oFTe)=\FTF\left(\bigcup \cMT \cap \oFTe \right) \subset 
\bigcup \FTF\left( \cMT \cap \oFTe \right),$$
\nid  the image of $\oFTe$ under 
$\FTF$ is a finite union of points belonging 
to $\cM$. Since, 
$\oFTe$ is invariant, this is a finite union of
points, and thus  a finite union of periodic
orbits with a finite period. The dynamics of neurons
$\in \bcF$ is driven by the periodic dynamics of firing neurons
and, from eq. (\ref{VjF}) it is easy to see that their trajectory converges to a constant. 
 \epr

\textbf{Remark.}
In the theorem, we have considered the case $\dAS=0$ as well. One sees
that there is no exponential proliferation of orbits after a finite time corresponding
to the time where all neurons satisfying property (1) 
have fired at least once. Indeed, then the reset
term project a convex domain onto a point, and this point cannot generate distinct orbits.
As discussed above the effect of $\cS$ is somehow cancelled by the reset intrinsic
to BMS model.
 Note however that there are at most $2^{NT}$  points in $\oM$,
and this number can be quite a bit large.

The situation is more complex if one cannot uniformly bound the first time of firing
as already discussed in section \ref{SGhost}.
 Assumptions (1), (2) of theorem  \ref{Tomega} leave us
on a safe ground but are they generic ? 
Let us now to consider the case where they are
not satisfied. Namely $\forall  \epsilon >0,  \forall T < \infty$,
$\exists \V \in \cM, \exists i \in \left\{1 \dots N \right\}$ such that 
$\forall t  \leq T, V_i(t) < \theta$ and $\forall t_0, \exists t \geq t_0$
such that $V_i(t) \geq \theta - \epsilon$. Call:

\beq\label{Bad}
\BTe=
\left\{
\V \in \cM | 
\exists i, \mbox{such that}:
\baR{ccc}
&(i)& \forall t \leq T, V_i(t) < \theta\\
&(ii)& \forall t_0, \exists t \geq t_0, V_i(t) \geq \theta- \epsilon.
\eaR
  \right\}
\eeq

We are looking for the set of parameters values $(\cW,\Ie)$ such that the 
set:

\beq\label{B}
\cB=\bigcap_{T > 0} \bigcap_{\epsilon >0} \BTe,
\eeq

\nid is non empty.  Note that $\BTep \subset \BTe$.
Thus, $\cB=\bigcap_{\epsilon >0} \Bie$.
We are thus looking for points $\V$ such that
$\forall t>0, V_i(t) < \theta$ and $\limsup_{t \to \infty} V(t)=\theta$.
Therefore, $\cB$ is exactly the set of ghost orbits.

We now prove that $\cB$ is generically empty.
Actually, we prove a more general result namely that $\dAS$ is generically non zero.
Before this, we have now to provide a definition of ``generic''.
For this, we shall assume from now on that the synaptic weights and inputs
belong to some compact space $\cH \subset \bbbr^{N^2+N}$.
 This basically means that the $W_{ij}$'s ($\Iei$'s)
are bounded (or have a vanishing probability to become infinite
if we deal with random matrices/inputs). 
One can endow $\cH$ with a probability measure
having a density with respect to the  Lebesgue measure.
This corresponds to choosing the synaptic weights 
and external currents with some probability distribution,
as we shall do in section \ref{SRS}.
We say that a subset  $\cN \subset \cH$ is ``non generic in a measure
theoretic sense'' if this set has zero measure.
This means that there is a zero probability to pick up a point 
in $\cN$ by choosing the synaptic weights and external currents randomly.
We say that it is ``non generic in a topological sense''
if it is the complementary set of a countable intersection
of dense sets \cite{KH}. This definition corresponds to the following
situation. If we find a point belonging to $\cN$ 
then a slight perturbation of this point leads out of $\cN$,
for any perturbation that belongs to an open dense set. 
In other words one can maybe find perturbations
that leave the point inside $\cN$ but they are specific and require e.g. precise algebraic
relations between the synaptic weights and/or input currents.
These two notion of genericity usually do not coincide \cite{KH}.

\bth\label{TGenA}
The subset of parameters $(\cW,\Ie) \in \cH$ such
that $\dAS=0$ is non generic in a topological and measure theoretic sense.
\enth

\textbf{Remark} Since this result holds for the two distinct
notions of genericity we shall use the term ``generic'' both 
in a topological and in a measure theoretic sense,
without further precision
in the sequel. 

\bpr

Take $\V \in \oM$ such that $d(\tV,\cS)=0$. 
Then, there exists
$i \in \sn $ such that $\inf_{t \geq 0} |V_i(t) - \theta|=0$.
We shall  consider separately two  cases.

\ben
\item Either $\exists B < \infty$  and a sequence  $\left\{t_k \right\}_{k \geq 0}$ 
such that $V_i(t_k)=\theta$ and $\delta_k<B, \forall k \geq 0$, where  
$\delta_k = t_{k+1}-t_k$.

\item Or $\V$ is a ghost orbit. This includes the case where $\delta_k$ defined above
is not bounded, corresponding to having $\lim_{t \to +\infty} V_i(t)=\theta$, but
also the case where $V_i(t)$ has no limit, and where $\limsup_{t \to +\infty} V_i(t)= \theta$ 
as in definition (\ref{DfGhost}).  
\een

\textbf{Case 1} According to  eq. (\ref{Vit}), the condition $V_i(t_{k+1})=\theta$ writes: 

\beq\label{Vsing}
V_i(t_{k+1})=
\sum_{n=0}^{\delta_k-1} \gamma^{n}
 (I^s_i(t_{k+1}-n-1)+\Iei)=\theta,
\eeq

\nid since $t_k$ is a firing time. Note that we have used
the notation $I^s_i(t)$ instead of the notation $I^s_i(\bEta_t)$,
used in eq. (\ref{Vit}), for simplicity.

The synaptic current $I^s_i$   takes only finitely many
values $\alpha_{i;l}=
\sum_{j \in \cD(\bEta_l)} W_{ij}$,
where $l$ is  an index enumerating the elements of $\cP$ ($l \leq 2^N$). 
Thus, the $\alpha_{i;l}$'s are \textit{only  functions of the $W_{ij}$'s} 
and they do not depend on the orbits. One can write:

\beq \label{Lindec}
\sum_{n=0}^{\delta_k-1} \gamma^{n}I^s_i(t_{k+1}-n-1)
=\sum_{l=1}^{2^N} \alpha_{i;l} x_{i;l}(t_{k+1}),
\eeq

\nid where:

\beq\label{xil}
x_{i;l}(t_{k+1}) = \sum_{n=0}^{\delta_k-1}
\gamma^n \chi\left[I^s_i(t_{k+1}-n-1) = \alpha_{i;l} \right],
\eeq

\nid where $\chi$ is the indicatrix function.
One may view the list $\left\{ x_{i;l}(t_{k+1})\right\}_{l=1}^{2^N}$ as the components of
a vector $\x_i(t_{k+1}) \in \bbbr^{2^N}$.
In this setting, relation (\ref{Vsing}) writes:

\beq \label{HPi}
\sum_{l=1}^{2^N} \alpha_{i;l} x_{i;l}(t_{k+1}) = \theta- \frac{1-\gamma^{\delta_k}}{1-\gamma}I_i^{ext},
\eeq

\nid since $I_i^{ext}$ does not depend on time. Equation (\ref{HPi}) defines
an affine hyperplane $P_{i,k}$ in $\bbbr^{2^N}$.

Call  $Q_{i,k}$ the set of  $x_{i;l}(t_{k+1})$'s. 
This is a finite, disconnected set, with  $\#Q_{i,k} = 2^{\delta_k}$
and whose elements are separated by a distance $\geq \gamma^{\delta_k}$.
Moreover, the $x_{i;l}(t_{k+1})$'s are positive. For each $k$ they obey:

\beq\label{Simplex}
\sum_{l=1}^{2^N} x_{i;l}(t_{k+1})=\sum_{n=0}^{\delta_k-1} \gamma^n = \frac{1 - \gamma^{\delta_k}}{1-\gamma}
\eeq 

This defines a simplex and  $Q_{i,k}$  belongs to this simplex.
Note  $Q_{i,k}$ \textit{does not} depend on the parameters $\cW,\Ie$.
However, the set of $x_{i;l}(t_{k+1})$'s values appearing
in eq. (\ref{HPi}) is in general a subset of $Q_{i,k}$ depending
on $(\cW,\Ie)$

Now, eq. (\ref{HPi}) has a solution if and only if $P_{i,k} \cap Q_{i,k} \neq \emptyset$.
Assume that we have found a point $R=(\cW,\Ie)$ in the parameters space $\cH$
such that  $P_{i,k} \cap Q_{i,k} \neq \emptyset$, for some $k$.
Since $Q_{i,k}$ is composed by finitely many isolated points, since the 
$\alpha_{i;l}$'s depend continuously on the $W_{ij}$'s and since
the affine constant of the hyperplane $P_{i,k}$ depends continuously of
$I_i^{ext}$, one can render the intersection  $P_{i,k} \cap Q_{i,k}$ empty
 by a generic (in both sense) small variation of the parameters $W_{ij}, I_i^{ext}$.
Therefore, the sets of points in $\cH$ such that $P_{i,k} \cap Q_{i,k} \neq \emptyset$, for some $k$,
is non generic.
Since we have assumed that the $\delta_k$'s are uniformly bounded by a constant
$B< \infty$, the condition $\exists k$ such that $V_i(t_k)=\theta$
corresponds to a finite union of non generic sets, and it is 
therefore non generic.  

Note that if $\delta_k$ is not bounded then the set of values $x_{i;l}=\sum_{n=0}^{\infty}
\gamma^n \chi\left[I^s_i(t_{k+1}-n-1) = \alpha_{i;l} \right]$ takes
uncountably many values. If $\gamma$ is sufficiently small this is Cantor set
and one can still use the same kind of argument as above.
On the other hand, if $\gamma$ is large this set fills
 continuously the simplex $\sum_{l=1}^{2^N} x_{i;l}= \frac{1}{1-\gamma}$
and one cannot directly use the argument above. More precisely one
must use  in addition some specificity of the BMS dynamics. This case is however
a sub case of ghost orbits. Therefore we  treat it
in the next item. \\

\textbf{Case 2.} We now prove that ghost orbits are non generic.
For this, we prove that if $R=(\cW,\I)$ is a point in $\cH$ 
such that the set $\cB$ defined by eq. (\ref{Bad}) is non empty, 
a small, generic, perturbation of $R$
leads to a point such that $\cB$ is empty. Thus, $\cB$ is generically empty
in both sense.

 Fix $\epsilon$ and take 
$\V \in \Bie$ (def. (\ref{Bad})). Then there is a $t_0$ such that
$\theta-\epsilon \leq V_i(t_0) < \theta$. Without loss of generality
(by changing the time origin) one may take $t_0=0$. Then, from eq. (\ref{Vit}),
$\forall t >0$,

$$\gamma^t(\theta-\epsilon)+\sum_{n=1}^t \gamma^{t-n}I_i(n-1) \leq V_i(t+1)
<\gamma^t \theta +\sum_{n=1}^t \gamma^{t-n}I_i(n-1),$$

\nid where we
have set $I_i(n-1) \equiv I_i(\eta_{n-1})$ to shorten the notations. Thus, $V_i(t)$ belongs
to an interval of diameter $\gamma^t\epsilon$.
Since $\epsilon$ can be arbitrarily small, and $t$ arbitrarily large we have
only to  consider the orbits such that $V_i(0)=\theta$, for some $i$.
There are finitely many such orbits.

Assume that $R=(\cW,\I)$ is such that $\cB$ is non empty.
Then, for some $i$, $\forall \epsilon>0$, there exists $t_0$ such that:

\beq \label{Cond1}
\theta-\epsilon \leq \gamma^{t_0} \theta +\sum_{n=1}^{t_0} \gamma^{{t_0}-n}I_i(n-1) <\theta,
\eeq
\nid and $\forall t >0$,
\beq\label{Cond2}
\sum_{n=1}^t \gamma^{t-n}I_i(n-1)<\theta(1-\gamma^t).
\eeq

Assume for the moment that there is only one neuron $i$ such that $\inf_{t \geq 0} |V_i(t)-\theta|=0$.
That is, all other neurons $j \neq i$ are such that $V_j(t)$ stays at a positive distance  
from $\theta$. In this case,
a small perturbation of the $W_{kj}$'s, where $k=1 \dots N$ \textit{but} $j \neq i$, or a small perturbation
of the $I_j^{ext}$'s \textit{will not}
change the values of the quantities $\eta_j(t)=Z(V_j(t))$, $t=0 \dots +\infty$. In this case, the current $I_i(n-1)$
in eq. (\ref{Cond1},\ref{Cond2}) does not change $\forall n \geq 0$.
Therefore there is a whole set of perturbations that do not remove the ghost orbit\footnote{For example,
there may exist submanifolds in $\cH$ corresponding to systems with ghost orbits.
A possible illustration of this is given in fig. \ref{Fdist}, section \ref{SRS}
where the sharp transition from a large distance $\dAS$
to very small distance $\dAS$ corresponds to a critical line
in the parameters space $\gamma,\sigma$ (see section \ref{SRS} for details). }.
But they are \textit{non generic} since a generic perturbation involves a variation
of \textit{all} synaptic weights $W_{kj}$ including $j=i$ and all currents
as well. 

Now, a small perturbation of some $W_{ki}$ or $I_i^{ext}$ has
the following effects. Call $V'_i(t)$ the perturbed value of
the membrane potential at time $t$.
\ben
\item Either $\forall t>0, V'_i(t) < \theta - \epsilon_0$, for some $\epsilon_0 > 0$. In this case,
condition (\ref{Cond1}) is violated and this perturbation has removed the ghost orbit.
Now, since $i$ is not firing, it does not act on the other neurons and we are done. 

\item Or there is some $t_0$ such that $V'_i(t_0) \geq \theta$. The condition (\ref{Cond1}) is violated and this perturbation 
also removes the ghost orbit. But, neuron $i$ is now firing and we have to consider
its effects on the other neurons. Note that the induced effects on neurons
$j \neq i$ is \textit{not small} since  neuron $j$ feels now,  at each times where
neuron $i$ fires, an additional term $W_{ji}$ which can be large.
Thus, in this case, a small perturbation induces drastic changes by ``avalanches'' effects.

Again, we have to consider two cases.

\ben

\item Either the new dynamical system resulting from this perturbation has no ghost
orbits and we are done.

\item Or, there is another neuron $i_1$ ($i_1 \neq i$) having a ghost orbit obeying conditions (\ref{Cond1},\ref{Cond2}).
But then one can remove this new ghost orbit by a new perturbation. Indeed,
as argued above, the fact that $i$ is now firing corresponds to
adding a term $W_{ji}$ to the synaptic current $I^s_j$ each time neuron
$i$ fires. Then, to still have a ghost orbit for $j$ one needs
specific algebraic relations between the synaptic weights and currents which
corresponds to a set of parameters of codimension lower than
$1$.  
The key point is that, following this argument, one can
find a family of generic perturbation that destroy the ghost orbits of $i_1$ without
creating again a ghost orbit for $i$. Then by a finite sequence of generic perturbations
one can find a point in $\cH$ such that $\cB$ is empty. 
\een
\een

Finally, we have to treat the case where more than one neuron are such that $\inf_{t \geq 0} |V_j(t)-\theta|=0$.
However these neurons correspond to case $1$ or to case $2$ and one can lead them to a positive distance from $\cS$
by a finite sequence of  generic perturbations.
\epr

%% file: SGenA.tex
\ssu{General structure of the asymptotic dynamics.} \label{SGenA}

We are now able to fully characterize the $\omega$ limit set
of $\cM$. 

\ben

\item\textbf{Neural death.} Assume that $\Iei < (1-\gamma)\theta$ and consider
the set  $\cMz= \left\{\V \ | V_i < \theta, \ \forall i \right\}$ corresponding
 to states where all neurons are quiescent. Under this  assumption on $\Iei$,
$\cMz$  is an absorbing domain ($\F(\cMz) \subset \cMz$) and
  $\F^t(\cMz) \to \frac{\Iei}{1-\gamma}$ as $t \to \infty$.
Thus, all neurons in this domain are in a ``neural death'' state in the sense that they never
fire. More generally, let $\cMe$ be a domain such that  $\exists t>0$ such that $\Ft(\cMe) \subset \cMz$
then all states  in $\cMe$ converge asymptotically to neural death (under the assumption $\Iei < (1-\gamma)\theta$).
Now, if  $\bigcup_{t \geq 0} \F^{-t}(\cMz) \supset \cM$ then all state $\forall \V \in \cM$ converges to neural death.
Such a condition is fulfilled if the total current is not sufficient to maintain a permanent neural activity.
This corresponds to the previous condition on $\Iei$ but also to a condition on the synaptic weights
$W_{ij}$. For example, an obvious, sufficient condition to have neural death is $\VM < \theta$.
More generally, we shall see in section \ref{SRS}, where random synapses are considered,
that there is a sharp transition from neural death to complex activity when the weights
have sufficiently large values (determined, in the example of section \ref{SRS}
by the variance of their probability distribution).

\item\textbf{Full activity.} On the opposite, 
consider now the domain $\cMu= \left\{\V \ | V_i \geq \theta, \ \forall i \right\}$ corresponding
 to states where all neurons are firing. Then, if $\forall i, \sum_{j=1}^N W_{ij} + \Iei \geq \theta$,
this domain is mapped into itself by $\F$ (where $\F(\cMu)$ is the point $\prod_{i=1}^N\sum_{j=1}^N W_{ij} + \Iei$)
 and all neuron fire at each time step, forever.
More generally, if $\bigcup_{t \geq 0} \F^{-t}(\cMu) \supset \cM$ then all state $\forall \V \in \cM$ converges to this state
of maximal activity. Such a condition is for example fulfilled if the total current is too strong.
\een 

These two situations are extremal cases that can be reached by tuning the total current. 
In between, the dynamics is quite a bit richer. One can actually distinguish $3$ typical situations 
described by the following theorem, which is a corollary of
 Prop. \ref{PFGen}, th. \ref{TAfini}, \ref{Tomega}
and previous examples.

\bth\label{TStructA}

Let 

\beq\label{Vsup}
\Vs = \max_{i=1 \dots N} \Vsi,
\eeq

\nid where:

\beq\label{Vsupi}
\Vsi= \sup_{\V \in \cM} \limsup_{t \to \infty}  V_i(t),
\eeq

\nid be the maximal membrane potential that the neurons can have
in the asymptotics. Then,

\ben
\item Either $\Vs < \theta$. Then $\Vs=\max_i \frac{\Iei}{1-\gamma}$, $\dAS=\theta-\Vs$ and 
$\oM$   is reduced
to a fixed point $\in \cMz$.  \textbf{[Neural death]}. 

\item Or $\dAS > \epsilon >0$ and $\Vs > \theta$. Then  $\oM$  is a finite
union of stable periodic orbits with a finite period \textbf{[Stable periodic regime.]}.

\item Or $\dAS=0$. Then necessarily $\Vs \geq \theta$.
 In this case the system exhibits a 
weak form of initial conditions sensitivity.  $\oM$ may contain ghost
 orbits but this case
is non generic. Generically, the $\omega$-limit set is a finite union
of periodic orbit.\textbf{[Unstable periodic regime.]}.
\een

\enth

\textbf{Remark} 

It results from these theorems that the BMS model is an automaton;
namely, the value of $\bEta$ at time $t$ can be written as a deterministic 
function of the past spiking sequences $\bEta(t-1), \bEta(t-2)$ etc ....
  However, the number of spiking patterns determining the actual value of $\bEta$
can be arbitrary large and even infinite, when $\dAS=0$. Moreover, the dynamics
is nevertheless far from being trivial, even in the simplest case $\gamma=0$
(see section \ref{SRS}).

%% file: SCoding.tex
\label{SCoding}

In this section we switch from the dynamics description in terms of orbit
to a description in terms of spiking patterns. For this
we first establish a relation between the values that the membrane potentials
have on $\oM$ and an infinite spiking patterns sequence, using
 the notion of \textit{global orbit} introduced in \cite{Bastien}.

%% file: SGlob.tex
\ssu{Global orbits.} \label{SGlob}

In (\ref{Vit}), we have implicitly fixed the initial time at $t=0$.
One can also fix it at $t=s$ then take the limit $s \to -\infty$. 
This allows us to remove the transients.  This leads to:

\beq\label{OGlob}
V_i(t)=\sum_{n=0}^{+\infty} \pi_i(n,t)\gamma^n I^s_i(t-n-1)
\eeq

\nid where:

\beq\label{pi}
\pi_i(n,t)=\prod_{k=0}^{n} 
\left(
1-\eta_{i;t-k-1} 
\right),
\eeq

\bdf
An orbit is \textit{global}  if there exists a legal sequence  
$\betat = \left\{\bEta_t\right\}_{t \in \bbbz} \in \SW$ such that $\forall t> 0$,
 $V_i(t)$ is given by (\ref{OGlob}).
\edf

\textbf{Remarks} 

\ben

\item In  (\ref{OGlob}) one considers  sequences
$\eta_{.;t-k-1}$ where $t-k-1$ can be negative, i.e.
$\left\{\eta_t\right\}_{t \in \bbbz} \in \SW$.
Thus a global orbit is such that its backward trajectory stays in $\cM$, $\forall t <0$. 

\item The quantity $\pi_i(n,t) \in \left\{0,1\right\}$, and
is equal to $1$ if and only if  neuron $i$, at time $t$,
has not fired since time time $t-n-1$.
Thus, if
$\tik $ is the last firing time, then $V_i(t)
=\sum_{n=0}^{t-\tik-1} \gamma^n I^s_i(t-n-1), \tik < t \leq \tkp $, is 
a a sum with a finite number of terms. The form (\ref{OGlob}) is
a series only when the neuron didn't fire in the (infinite)
past.

\een

Denote by $\cG$ the set of global orbits.
The next theorem is an (almost) direct transposition of proposition 5.2 proved
by Countinho et al. in \cite{Bastien}. 
 However, the paper \cite{Bastien}
deals with a different model and slight adaptations of the proof have to be made.  
The main difference is the fact that, contrarily to their model, it is not true
that every point in $\bbbr^N$ has a uniformly bounded number of pre-images.
This is because $\F$ typically project a domain onto a domain of lower dimension
in all directions where a neuron fires (and this effect is not equivalent to setting
$a=0$ in  \cite{Bastien}). Therefore, to apply    Countinho et al. proof
we have to exclude the case where a point has infinitely many pre-images.
But it is easy to see that in the generic situation of th. \ref{Tomega} 
 any point of  $\oM$ has a finite number of pre-images in $\oM$ (since $\oM$
has finitely many points).

The version of  Countinho et al. theorem for the BMS model is therefore.

\bth \label{TGlob}
.
$\oM=\cG$ for a generic set of $(\cW,\Ie)$ values.
\enth

\textbf{Remark} For technical reasons we shall consider the attractor $\cA$ definition (eq. \ref{A})
 instead of the $\omega$-limit set. But these two notions coincide
whenever there is no ghost orbit (generic case).

\bpr
The inclusion $\cG \subset \cA$ is proved as follows.
Let   $\V \in  \cG$ and $\tV=\left\{\V(t) \right\}_{t \in \bbbz}$
be the corresponding global orbit. 
Since, $\forall t, \ n$,
$$\min_{i} \sum_{j=1}^N W_{ij} \leq
I_i^s(t-n-1) 
\leq \max_{i} \sum_{j=1}^N W_{ij},$$
 one has
$$\sum_{n=0}^\infty \gamma^n \left(\min_{i} \sum_{j=1}^N W_{ij}+\Iei\right) \leq 
V_i(t) 
\leq \sum_{n=0}^\infty \gamma^n \left(\max_i \sum_{j=1}^N  W_{ij} +\Iei\right)$$
$$\Rightarrow
\Vm = \leq 
V_i(t) 
\leq 
\VM.$$
Therefore, 
$\V(t) \in \cM \subset \bmd, \ \forall t \leq 0, \ \delta > 0$.
Hence $\V \in \bigcap_{t=0}^\infty
\Ft(\bmd)$ and $\cG \subset \bigcap_{t=0}^\infty
\Ft(\bmd)$. From (\ref{ContBoule}),
$\bigcap_{t=0}^\infty
\Ft(\bmd) \subset \cA$, and $\cG \subset \cA$.\\

The reverse inclusion $\cA \subset \cG$ is a direct consequence of
the fact that any point of $\cA$ has a pre-image in $\cA$. 
Therefore, $\forall \V \in \cA$, one can construct
an orbit $\left\{\V(t) \right\}_{t \leq 0}$ such that
$\V(0)=\V$, $\V(t+1)=\F(\V(t))$ and $\V(t) \in  \cA, \
\forall t \leq 1$. This (backward) orbit belong to $\cM$ and
the value of  $\V(t)$ is given by (\ref{OGlob}).
Thus $\V \in \cG$, so $\cA \subset \cG$.
\epr

\textbf{Remark.} Theorem \ref{TGlob} states that each  point in the attractor is generically encoded
by a legal sequence $\betat$. This is one of the key results of this contribution.
Indeed, as discussed in the introduction, the ``physical'' or ``natural'' quantity 
for the neural network is the membrane potential. However, it
is also admitted in the neural network community that the information
transported by the neurons dynamics is contained in the sequence of spikes
emitted by each neurons.  In the BMS model
such a sequence is exactly given by  $\betat$ since
on the $i$-th line $\eta_{t;i}$ one can read  the sequence of spikes
(and the firing times) emitted by  $i$. The theorem establishes 
that, in the BMS model, it is equivalent to consider the membrane potentials
or the spiking sequences: the correspondence is one to one.
This suggests a ``change of paradigm'' where one switches from
the dynamics of membrane potential (eq. \ref{DNN}) to the dynamics
of spiking patterns sequences. This is the point of view developed
in this series of papers, where some important consequences 
are inferred.

%% file: SRS.tex
\ssu{Random synapses.} \label{SRS}

In this paper we have established general results on the BMS model dynamics, and we have established
theorems holding either for all possible values of the $W_{ij}$'s and $\Iei$'s or for a generic set.
However, and obviously, the dynamics exhibited by the system (\ref{DNN}) depend on
the matrix $\cW$ (and the input $\Ie$) and quantities such as $\dAS$ or $\Vs$ in th. \ref{TStructA}
are dependent on these parameters. A continuous variation of some $W_{ij}$ or some $\Iei$
will induce quantitative changes in the dynamics (for example it will reduce the period or
the number or periodic orbits). It is therefore interesting to figure out what are the regions
in the parameters spaces $\cW,\Ie$ where the dynamics exhibits a different quantitative behaviour.

A possible way to explore this aspect is to choose  $\cW$ (and/or $\Ie$) randomly, with some
probability $\PrW$ ($\PrI$) having a density. A natural  starting point is the use of Gaussian independent, identically
distributed variables, where one
varies the statistical parameters (mean and variance). Doing these variations,
one performs sort of a fuzzy sampling of the parameters space, and one somehow expects
the behaviour observed for a given value of the statistical parameters to be characteristic
of the region of $\cW,\Ie$ that the probabilities $\PrW,\PrI$ weight (more precisely, one expects to observe
a ``prevalent'' behaviour in the sense of Hunt \& al. \cite{Sauer}). 

Imposing such a probability distribution  has several consequences.
First, the synaptic currents and the membrane potentials become random variables
whose law is induced by the distribution $\PrWi=\PrW\PrI$ 
and this law can be somehow determined \cite{BMSRand}. But, this has another,
more subtle effect.
Consider the set $\SL$ of \underline{all} possible  sequences on $\Lambda=\sn$. Among them,
the dynamics (\ref{DNN}) selects a subset of legal sequences, $\SW$, defined by the compatibility
conditions (\ref{Transitions}) and the transition graph $\GW$. Thus, changing $\cW$ ($\Ie$) has the effect
of changing the set of legal transitions that the dynamics selects. 
From a  practical point of view, this simply means that the typical raster plots observed in the asymptotic
dynamics depend on the $W_{ij}$'s and on the external current $\Ie$. This remark is somewhat evident. 
However, a question is \textit{how} the statistical parameters of the distribution $\PrWi$ 
acts on the dynamics typically observed in the asymptotics (e.g. how it acts on the parameters
$\Vs,\dAS$). This question can be addressed by combining the dynamical system approach of
the present paper, probabilistic methods and mean-field approaches
from statistical physics (see \cite{CS,SC} for an example of such
combination applied to neural networks).  A detailed description of this aspect would increase consequently
the size of the paper, so this will be developed in a separate work \cite{BMSRand}.
Instead, we would like to briefly comment results obtained by BMS.\\

Indeed, the influence of the statistical parameters of the probability
distribution of synapses on the dynamics has been investigated by BMS, using a different approach
than ours.
They have considered the case where the $W_{ij}$'s are Gaussian with
zero mean and a variance $\sigma^2$, and where the external current was zero.
 By using a mean-field approach they were able to
obtain analytically a (non rigorous) self-consistent
 equation (mean-field equation) for the probability $x_t$ that a neuron fires
at a given time. This equation always exhibits the locally stable solution $x=0$ corresponding to the ``neural death''.
For sufficiently large $\sigma$ another stable solution appears by a saddle-node bifurcation,
corresponding to a non zero probability of firing. In this case, one has two stable
coexisting regimes (neural death and non zero probability of firing), and one reaches
one regime or the other according to the initial probability of firing. Basically,
if the initial level of firing is high enough, the network is able to maintain
a regime with a neuronal activity. This situation appears for a sufficiently large
value of $\sigma$, corresponding to a critical line in the plane $\gamma,\sigma$.
The analytical form of this critical line was not given by BMS. Moreover, the
mean-field approach gives information about the average behavior of an ensemble of neural networks
in the limit $N \to \infty$. The convergence involved in this limit is weak convergence (instead of almost-sure
convergence). Therefore, it does not tell us what will be the typical behaviour of one infinite sized neural network.
Finally, the mean-field approach does not allow to describe 
the typical dynamics of a finite sized network.

To study the finite size dynamics BMS used numerics and gave evidence
of three regimes. 

\bit
\item \textbf{Neural death.} After a finite time the neurons stop to fire.
 
\item \textbf{Periodic regime.} This regime occurs when $\sigma$ is large enough.

\item \textbf{``Chaos''.} Moreover, BMS exhibit an intermediate regime, between
neural death and periodic regime,
that they associate to a chaotic activity. In particular,
 numerical computations with the Eckmann-Ruelle algorithm \cite{ER}
exhibit a positive Lyapunov exponent. This exponent
decreases to zero when $\sigma$ increases, and becomes negative in the periodic
regime.
\eit

\begin{figure}[ht]
\begin{center}
\includegraphics[height=8cm,width=8cm,clip=false]{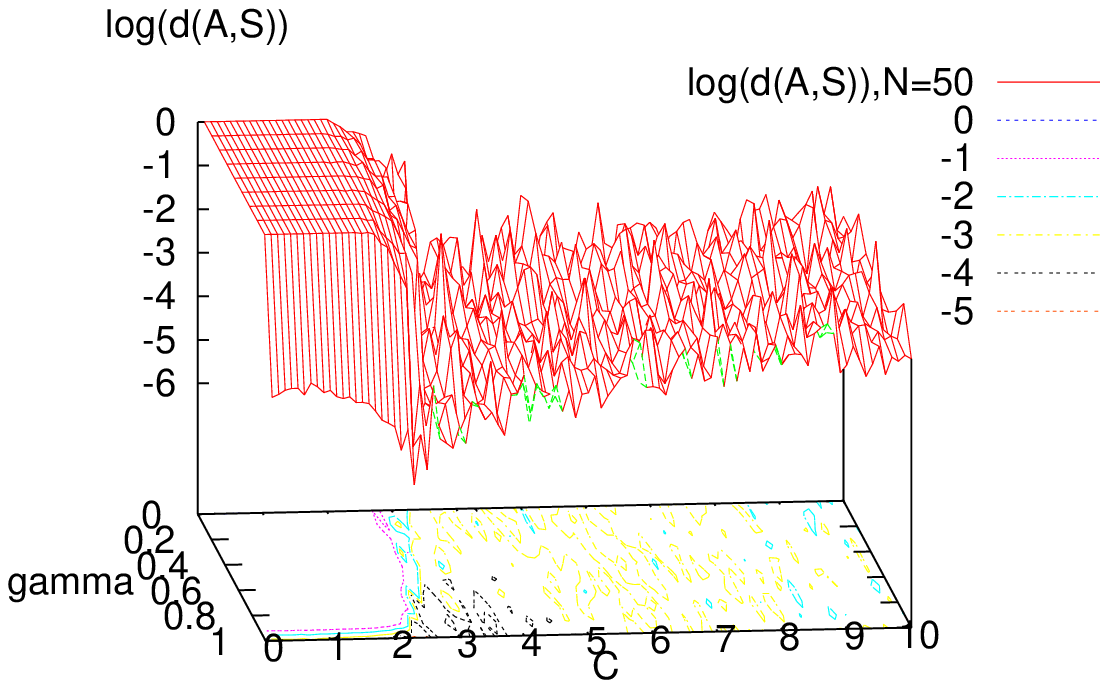}
\includegraphics[height=8cm,width=8cm,clip=false]{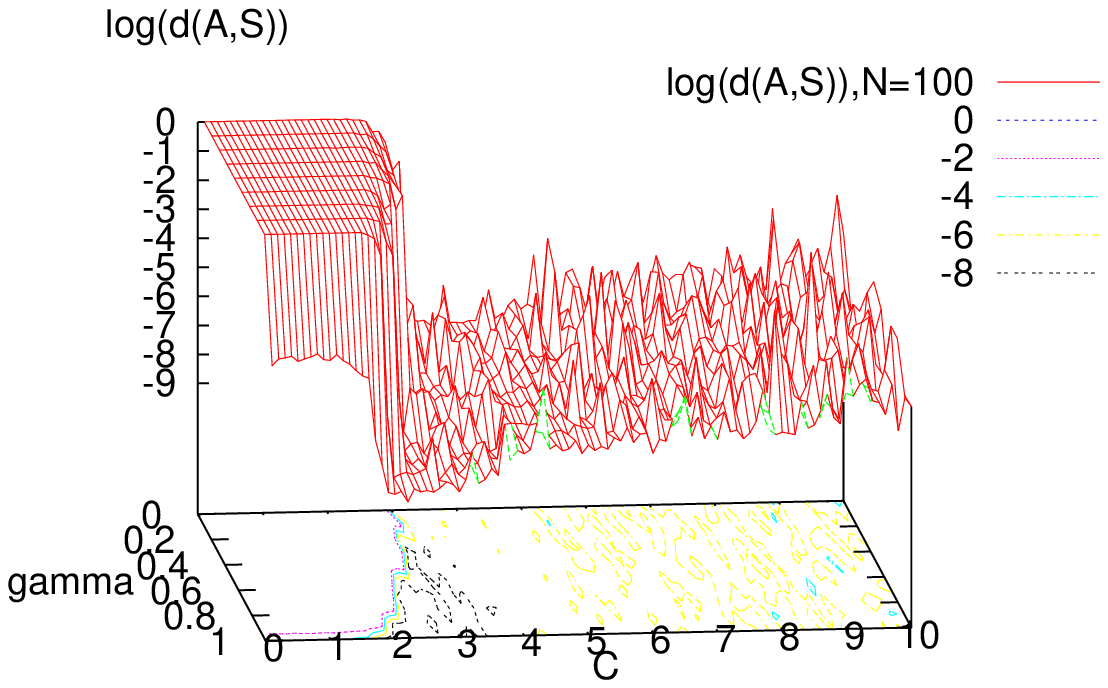}
\vspace{0.5cm}
\caption{Fig. \ref{Fdist}a. Average value of the distance $d(\cA,\cS)$ versus $\gamma,C$, for $N=50$.
Fig. \ref{Fdist}b. $N=100$ (in $\log_{10}$ scale).}
\label{Fdist}\end{center}
\end{figure}

Their conclusion concerning the existence of a chaotic regime is in contradiction with
theorem \ref{TStructA}. We would like now to briefly comment this contradiction
(a more detailed investigation will be done in \cite{BMSRand}). The fig. \ref{Fdist}a,b
presents the results of a numerical simulation computing
the average distance $\dAS$ as a function of $\gamma$ and of the variance of the synaptic weights.
More precisely, we have considered, as BMS, the case of Gaussian independent, identically
distributed random $J_{ij}$'s, with zero expectation and variance $\sigma^2=\frac{C^2}{N}$. 
(We have adopted the standard scaling of the variance with $\frac{1}{N}$. Indeed, in the present
case the neural network is almost surely fully connected and the
scaling $\frac{C^2}{N}$ is used in order that the probability of
the total currents $I_i$  has a  variance independent of $N$). 

Clearly, the average distance becomes very small  when $C$ crosses
a critical line in the plane $C,\gamma$. 
However, in the numerical
experiments of Fig. \ref{Fdist} the smaller measured value for the distance is $\sim 10^{-8}$
for Fig. \ref{Fdist}b, corresponding to a very large characteristic time
 well beyond the transients usually considered in the numerics (eq. (\ref{Td}).
Moreover,  the average distance approaches zero rapidly as $N$ growths. Thus, there is
sharp transition from neural death to chaotic activity in the limit $N \to \infty$,
when crossing a critical line in the plane $C,\gamma$ (``edge of chaos''). This line
can be determined by mean-field methods analogous to those used in \cite{JP} and
corresponds to the transition found by BMS \cite{BMSRand}. 
In fig. \ref{Fdist}a,b, one also remarks that after the transition
 $\dAS$ growths slowly when $C$ increases. For the illustration
of this aspect we have drawn the log of the distance in fig. \ref{Fdist}a,b.

Hence, for finite size $N$ the situation is the following. Start from a small variance
parameter $C$ and increase it, and consider the stationary regime
typically observed. There is first a neural death regime. After this,
 there is a regime where the dynamics has a large number
of periodic orbits and very long transients.
 This regime is numerically indistinguishable from chaos\footnote
{
 Moreover, it is likely
that the phase space structure has some analogies with spin-glasses\cite{BMSRand}.
 For example, if $\gamma=0$ the dynamics is essentially equivalent to the Kauffman's cellular
automaton \cite{Kauff}. It has been shown by Derrida and coworkers \cite{DF},\cite{DP}
that the Kauffman's model has a structure similar to the Sherrington-Kirckpatrick spin-glass
model\cite{MPV,Sherrington}. The situation is even more complex when $\gamma \neq 0$.
It is likely that we have in fact a situation very similar to discrete time neural networks
with firing rates where a similar analogy has been exhibited \cite{EPL},\cite{JP}.
}.
In particular, usual numerical methods, computing Lyapunov exponents by studying the behaviour
of a small ball of perturbed trajectories centered around a mother trajectory, will
find a positive exponent. Indeed,  if the size $\eta$ of this ball is larger than the distance
$\dAS$ one will observe an effective expansion and initial condition
sensitivity, as argued in the section \ref{SSing}. This will result
in the measurement of an effective positive Lyapunov exponent, stable
with respect to small variation of $\eta$, as long as $\eta>>\dAS$.
Though this exponent is, strictly speaking, spurious, it captures
the most salient feature of the model: sensitivity to perturbations
with a finite amplitude.
When $C$ increases further, the distance to the singularity set
increases. There is then a $C$ such that the typical periodic orbit
length becomes of the order of magnitude of the time range used
in the numerical simulation, and one is able to see that dynamics is periodic.\\

In the light of this analysis we claim that BMS results are essentially correct 
though we have shown
that there is no strictly speaking chaotic regime. Moreover, they are, in some sense, more
relevant than theorems
 \ref{TGenA},\ref{TStructA} as far as numerics and practical aspects are concerned. 
However, the analysis of the present paper permits to have a detailed description of the typical
dynamics of a given finite sized network (without averaging), based on rigorous results.
This is useful when dealing with synaptic plasticity and learning effects where a \textit{given} pattern
is learned in a \textit{given} network. (This aspect is shortly discussed below and will be developed elsewhere).

%
%

%% file: SNoise.tex
\ssu{Adding noise to the dynamics.}\label{SNoise}

It is usual in neural network modeling to add Brownian noise to the
deterministic dynamics. This noise accounts for different effects
such as the diffusion of neurotransmitters involved in the synaptic
transmission, the degrees of freedom neglected by the model,
external perturbations, etc ... Though it is not evident that the ``real  noise''
is Brownian, using this kind of perturbations has the advantage
of providing a tractable model where  standard theorems in the theory
of stochastic processes \cite{FauTou} or methods in non equilibrium statistical physics (e.g. Fokker-Planck
equations \cite{BrunHak}) can be applied.

The addition of this  type of noise to the dynamics of BMS model will
result, in the region where $\dAS$ is small, in an effective initial
condition sensitivity and an effective positive Lyapunov exponent.

More precisely, consider a noisy version of (\ref{DNN}).

\beq\label{BMSNoise}
V_i(t+1)=\gamma V_i(t) \left(1 - Z[V_i(t)] \right)+ \sum_{j=1}^N W_{ij}Z[V_j(t)]+ \Iei(t) + B_{i}(t); \qquad i=1 \dots N.
\eeq
 
\nid where $\B \deq \left\{ B_{i}(t)\right\}_{i=1;t=0}^{N,\infty}$ is a Gaussian random
process with zero mean and  a covariance $Cov(B_i(t),B_j(s)) = \sB \delta_{t,s}\delta_{i,j}$.
The probability distribution of the stochastic process $\V$,
on a finite time horizon $T$, for a fixed realisation of the $\cW$ can be obtained
by using a discrete time version of Girsanov theorem \cite{Skor},\cite{SC}.
From this, it is possible to estimate the probability that a trajectory approaches
the singularity set $\cS$ within a finite time $T$ and a distance $d$ by using
Freidlin-Wentsel estimates \cite{FW}. Also, eq. (\ref{CBC}) is useful to estimate the
measure of points having a local stable manifold. In this context one can compute the probability to approach the
singularity set within a distance $\epsilon$; also one can construct a Markov chain for the transition
between the attraction basin of the periodic orbits of the unperturbed dynamics. 
This will be done in a forthcoming paper.

%% file: SIt.tex
\ssu{Time dependent input.}\label{SIt}

One may also wonder what happens to the present analysis when a deterministic,
time dependent external input, is imposed upon the dynamics (the case
of a stochastic input is covered by eq. (\ref{BMSNoise}) above).  
Away from the singularity set ($\dAS$ large) the effect of a time dependent
input with a small amplitude (lower than $\dAS$) will not be different from the
case studied in the present paper. This is basically because a small
input may be viewed as a perturbation of the trajectory, and the contraction
properties of the dynamics will damp the perturbation as long as
the trajectory stays away from the singularity set. 

The situation is different if, at some place, the action of the time dependent
input leads to a crossing of the singularity. This crossing can basically occur
with a time independent input, but in the time dependent case there is a particularly
salient effect, that may be easily revealed with periodic external currents.
That is \textit{resonance effects}. If the unperturbed trajectory has some typical
recurrent time to come close to the singularity set, and if the time dependent perturbation
is not synchronized with this recurrence time, one expects that the contraction effect
will damp the perturbation with no clear cut ``emergent'' effect. On the other hand, if the period
of the periodic signal is a multiple of the recurrence time, there may be a major
effect. The result would be a frequency dependent response of the system exhibiting
sharp peaks (resonance). This statement is actually more than a conjecture. Such resonances
effects have indeed been exhibited in a recurrent discrete time neural network with firing
rates \cite{JAB1},\cite{JAB2},\cite{JAB3}. It has been shown that applying a periodic input
is a way to handle the interwoven effects of non linear dynamics and synaptic topology.
Similar effects should be observed in BMS model.

%% file: SHebb.tex
\ssu{Learning and synaptic plasticity.} \label{SHebb}

What would be the effect of a synaptic weight variations (synaptic plasticity, LTD, LTP, STDP,
Hebbian learning) on the dynamical system (\ref{DNN}) ? 
These variations corresponds
to moving the point  corresponding to the dynamical system
in the parameters space $(\cW,\Ie)$. This motion is neither random nor
arbitrary. Indeed, assume that one imposes to the neural network
an input/stimulus $\Ie=\left\{\Iei(t)\right\}$. $\Iei$ modifies directly the level of activity 
of  neuron $i$, and acts indirectly on other neurons
(provided that the synaptic graph is connected). A simple stimulus can therefore
strongly modify the dynamics, the attracting set, the distance $\dAS$, etc .... 

In the case where $\Ie$ does not depend on time, the following result follows directly from
the analysis presented in this paper.

\bth\label{ThOrbInp} 
For a generic set of values of $(\cW,\Ie)$,
there exists a finite partition of  $\cM=
\bigcup \cD_n$, such that $\forall \V \in \cD_n$ the $\omega$-limit
set of $\V$, $\omega(\V)$ is a stable periodic orbit,
with a finite period. This orbit depends on $\Ie$.
\enth 

\bpr
$\oM$ is generically a finite union of periodic orbits
with a finite period.
Each of orbit $n$ has an attraction basin $\cD_n$
and the attraction basins consitute a partition of 
$\cM$. 
\epr

This orbit (resp. its coding) may be viewed as the
dynamical response of the neural network to input $\Ie$, 
whenever the initial conditions are chosen $\in \cD_n$.
In this way, the neural network associates to an input
a dynamical pattern encoded in the spiking sequence of this periodic orbit. In the same way one 
can associate to a series of inputs a series of  periodic orbits (resp. codes),
each orbit being specifically related to an input. This property
results directly from th. (\ref{ThOrbInp}) without particular assumption on the $W_{ij}$'s.\\

However, there might exist a large number of domains $\cD_n$ and a large number
of possible responses (orbits). Moreover, an orbit can be complex, with a very long period.
This is particularly true at the ``edge of chaos''. Indeed, consider the case where 
the distance $\dAS$ is small, when the input is present. Then, dynamics is indistinguishable from chaos and
the dynamical ``signature'' of the input is a very complex orbit, requiring a very long time to
be identified. In other words, if one imagine a layered structure where the present  neural network acts as a retina
and where another neural network is intended to identify the orbit and 
``recognize'' the input, the integration time of the retina will be very long at the edge
of chaos. On the opposite,  one may expect that a learning phase allows this system
to associate the input to an orbit with a simple structure (small period) allowing a
fast identification of the input. 

It has been shown, in the case of recurrent neural networks with a sigmoidal transfer function \cite{Dauce}, that Hebbian learning leads to a reduction of
chaos towards a less complex dynamics, permitting to associate a pattern to 
simple orbits. The same effect has been observed by BMS \cite{BMS} applying an STDP like rule
to the model (\ref{DNN}). In both cases, it has been observed that a synaptic
evolution (Hebb or STDP) leads to associate to the input a sequence of
orbits whose complexity decreases during the synaptic weights evolution.
In the present context, this suggests that $d(\omega(\V),\cS)$ increases during
this evolution 
(note that the  evolution is entirely dependent on the initial  condition, $\V$. ).\\

A related question is: how do the statistical properties of raster plots evolve during synaptic weights
evolution ? This question, and more generally the effect of synaptic evolution on dynamics
can be addressed using tools from dynamical systems theory, in the spirit of the present paper.
This will be the subject of a forthcoming paper. However, in the next section we mention
briefly how tools from ergodic theory (thermodynamic formalism) can be used.

%% file: SGibbs.tex
\ssu{Statistical properties of  orbits.}\label{SStat}

As we saw, the dynamics of (\ref{DNN}) is a rather complex and can be,
from an experimental point of view, indistinguishable from chaos.
Consequently, the study of the finite evolution of the membrane
potential (resp. the spiking patterns sequence) does not tell us what
will be the further evolution, whenever the time of observation
is smaller than the characteristic time $T_d$ of eq. (\ref{Td}). In this
sense, the system is producing entropy on a finite time horizon. Thus,
provided that $\dAS$ is sufficiently small, one can do ``as if'' the system
were chaotic and use the tools for analysis of chaotic systems. This also
holds when one adds noise on the dynamics. A particularly
useful set of tools is provided by ergodic theory and the thermodynamic
formalism. In this approach one is interested in the statistical behavior
of orbits, characterized by a suitable set of probability measure.
A natural choice are Gibbs measures in the sense of Sinai-Ruelle-Bowen
\cite{SRB}. In a forthcoming paper we indeed show that Gibbs
measures arise naturally in BMS model.  They come either from statistical
inference principles where one tries to maximize the statistical
entropy given a set of fixed given quantities such as correlations functions
or mean firing rate (a prominent example of application of this principle
is given in \cite{Schneid}). They also arise when one wants to
study the effect of synaptic plasticity (learning, STDP) on the selection
of orbits. In the context of BMS model one can show
that Hebbian learning and STDP are related to a variational principle
on the topological pressure, which is the analogon of free energy
in statistical mechanics.  

%% file: Sdt.tex
\ssu{The limit $dt \to 0$.}\label{Sdt}

In the definition of the BMS model, one uses a somewhat rough approximation
consisting in approximating the differential equation of the Integrate and Fire
model with a Euler scheme, and discretizing time. A central question is:
what did we lose by doing this, and is the model still relevant as a
neural network model ? As mentioned in the introduction, this requires developments
done elsewhere \cite{CV}. But we would like here to point 
out here a few remarks on this aspect.

\bit

\item From the ``biological'' point of view 
the Integrate and fire model with continuous time is already a rough approximation
where the characteristic time for the neuron response is set to zero. One can actually
distinguish (at least) $3$ characteristic time scales in  neuron dynamics descriptions
based on differential equations. The ``microscopic time'' $dt$ corresponds somehow
to the shortest time scale involved in the spike generation (e.g. microscopic mechanism
leading to opening of ionic channels). The ``reaction time'' $\tau_r$ of the neuron corresponds
to the time of raise and fall for the spike. If one focuses on spikes (and does not consider
time averaging over sliding windows leading to the firing rate description)
the last relevant time scale is the characteristic time $T$ 
required  for the neural network to reach a stationary regime. One expects to
have $dt << \tau_r << T$. In the IF model, however, the time reaction $\tau_r$
is considered to be instantaneous (thus $\tau_r \leq dt$). This leads to delicate problems for the definition
of the time of firing and requires the introduction of the ``$t^-$ notation''.
Using a discrete time approximation allows to circumvent this problem and corresponds
somehow to pose $dt=\tau_r=1$. 

One may reject this procedure a priori. Our philosophy is instead to extract as much results
as possible from the discrete time spiking model and decide a posteriori what has been
lost (or won).

\item  From the dynamical system point of view, the limit $dt \to 0$ raises two 
problems. On one hand, the trajectories become continuous. Then one may have situations
where the trajectory accumulates on $\cS$ and where a small variation of the $W_{ij}$'s
is not able to remove the intersection (as it is the case in th. \ref{TGenA}). This type
of situation is known in the field of genetic networks (see \cite{Farcot} and references therein).
However, as mentioned in the paper, the situation is slightly different here, because
of the neurons reset, leading to an infinite contraction of a domain onto a point. 
This effect really simplifies the dynamics study, and is still present in the continuous
time case. However, this aspect would require careful investigations, not in the scope
of the present work. 

The second problem is the use of a Euler scheme in the discretization. Using
more elaborated schemes would complicate the analysis since the model
would loose its convenient piecewise affine structure. We don't know
what this would add.

\item Finally, from a numerical point of view, softwares use discrete time.
One aspect that interests us particularly is to know what are actually
the computing capacities of the discrete time model compared to classical
IF models and how much has been lost.
  
\eit

%% file: Ack.tex
\begin{acknowledgements}
I would like to thank G. Beslon, O. Mazet, H. Soula and M. Samuelides who told
me about the ``BMS'' model. I especially thank H. Soula for fruitful exchanges and M. Samuelides, J. Touboul and E. Ugalde
for a careful reading of this work.
This paper greatly benefited from intensive discussions with B. Fernandez, A. Meyronic,
and R. Lima. The remarks, questions and suggestions of T. Vi\'eville were
determinant in the redaction of this paper. I warmly acknowledge him.
Finally, I am grateful to the referees for helpful remarks and constructive
criticism.
 
\end{acknowledgements}